\documentclass[11pt,reqno]{amsart}

\usepackage{amssymb,amsmath,eufrak}
\author[Florent Benaych-Georges]{Florent Benaych-Georges}
\title[Rectangular random matrices, entropy, and Fisher's information]{Rectangular random matrices, related free entropy and free Fisher's information}
\date{\today}
\newcommand{\Alg}{\operatorname{Alg}}
\newcommand{\Haar}{\operatorname{Haar}}
\newcommand{\Proj}{\operatorname{Pr}}
\newcommand{\diag}{\operatorname{diag}}
\newcommand{\Adj}{\operatorname{Adj}}

\newcommand{\NC}{\operatorname{NC}}

\newcommand{\ED}{\operatorname{E}}
\newcommand{\EDn}{\operatorname{E_{\mathit{n}}}}
\newcommand{\EDpi}{\operatorname{E_{\pi}}}
\newcommand{\EDsigma}{\operatorname{E_{\sigma}}}

\newcommand{\ds}{\displaystyle}

\newcommand{\tr}{\operatorname{tr}}

\newcommand{\Span}{\operatorname{Span}}
\newcommand{\ELinf}{\operatorname{E_{\mathit{L}^{\infty}}}}
\newcommand{\mf}{\mathfrak}
\newcommand{\Mn}{\mathfrak{M}_n}

\newcommand{\Tr}{\operatorname{Tr}}
\newcommand{\ninf}{\underset{n\to\infty}{\longrightarrow}}
\newcommand{\ssi}{if and only if }

\newcommand{\teo}{theorem }

\newcommand{\E}{\mathbb{E}}

\newcommand{\R}{\mathbb{R}}
\newcommand{\C}{\mathbb{C}}
\newcommand{\n}{\mathbb{N}}

\newcommand{\ud}{\mathrm{d}}

\newcommand{\pro}{probability }

\newcommand{\f}{\frac}
\newcommand{\ff}{\frac{1}}
\newcommand{\lf}{\left}
\newcommand{\ri}{\right}

\newcommand{\st}{such that }
\newcommand{\la}{\lambda}
\newcommand{\La}{\Lambda}
\newcommand{\tii}{\scriptstyle\times\ds\!}

\newcommand{\vfi}{\varphi}
\newcommand{\ste}{\, ;\, }
\newcommand{\mc}{\mathcal }
\newcommand{\eps}{\varepsilon}
\newcommand{\arc}{{\scriptstyle\boxplus_{\la}}\!}

\newcommand{\stt}{\scriptstyle\times\ds\!}
\newcommand{\A}{\mc{A}}
\newcommand{\D}{\mc{D}}
\newcommand{\B}{\mc{B}}

\newcommand{\fao}{free with amalgamation over $\D$ }
\newcommand{\faop}{free with amalgamation over $\D$. }
\newcommand{\faose}{free with amalgamation over $\D$}

\numberwithin{equation}{section}
\newtheorem{Th}{Theorem}[section]
\newtheorem{propo}[Th]{Proposition} 
 
\newtheorem{lem}[Th]{Lemma}

\newtheorem{rmq}[Th]{Remark}
\newtheorem{cor}[Th]{Corollary}
\newtheorem{Def}[Th]{Definition}

\newenvironment{pr}{\noindent {\bf Proof. }}{\ \ \ $\square$}

\long\def\symbolfootnote[#1]#2{\begingroup
\def\thefootnote{\fnsymbol{footnote}}\footnote[#1]{#2}\endgroup}

\begin{document}
\maketitle
\symbolfootnote[0]{{\it MSC 2000 subject classifications.}  primary 
46L54, 
secondary 15A52} 
 
\symbolfootnote[0]{{\it Key words.} free probability, random matrices, free entropy, free Fisher's information, Marchenko-Pastur distribution}

\begin{abstract}We prove that independent rectangular random matrices, when embedded in a space of larger square matrices, are asymptotically free with amalgamation over a commutative finite dimensional subalgebra  $\D$ (under an hypothesis of unitary invariance). Then we consider elements of a finite von Neumann algebra containing $\D$, which have kernel and range projection in $\D$. We associate them a free entropy with the microstates approach, and a free Fisher's information with the conjugate variables approach. Both give rise to optimization problems whose solutions involve freeness with amalgamation over $\D$. It could be a first proposition for the study of operators between different Hilbert spaces with the tools of free probability. As an application, we prove a result of freeness with amalgamation between the two parts of the polar decomposition of $R$-diagonal elements with non trivial kernel. 
\end{abstract}
\tableofcontents

\section*{Introduction} In a previous paper (\cite{fbg.AOP.rect}), we considered an independent  family of rectangular random matrices with different sizes, say $n\tii p$, $p\tii n$, $n\tii n$ and $p\tii p$. We embedded them, as blocks, in $(n+p)\tii (n+p)$ matrices by the following rules\begin{equation}\label{je.chancelle.un.tt.petit.peu.quand.je.te.rencontre} M\to \ds \lf\{
 \begin{array}{clcl}\begin{bmatrix}M&0\\ 0&0\end{bmatrix}&\textrm{ if $M$ is $n\tii n$,}&\begin{bmatrix}0&M\\ 0&0\end{bmatrix}&\textrm{ if $M$ is $n\tii p$,}\\ &&&\\
\begin{bmatrix}0&0\\ 0&M\end{bmatrix}&\textrm{ if $M$ is $p\tii p$,}&
\begin{bmatrix}0&0\\ M&0\end{bmatrix}&\textrm{ if $M$ is $p\tii n$,}\\
\end{array}\ri.\end{equation}  and we proved that under an assumption 
of invariance under actions of unitary groups and of convergence of singular laws (i.e. uniform distribution on eigenvalues of the absolute value), the embedded matrices are asymptotically free with amalgamation on the two-dimensional commutative subalgebra $\D$ generated by the projectors $$\begin{bmatrix}I_n&0\\ 0&0\end{bmatrix} ,\quad \begin{bmatrix}0&0\\ 0&I_p\end{bmatrix}.$$ Asymptotically refers to the limit when $n,p\to\infty$ in a ratio having a non negative limit. In fact, we considered not only two sizes $n, p$, but a finite family $q_1(n),\ldots, q_d(n)$ of sizes, and the large matrices where represented as $d\tii d$ block matrices. 

In this paper, we prove a similar result with different technics (which allows us to remove the hypothesis   of convergence of singular laws). Then we consider a $W^*$-\pro space $(\A,\vfi)$ endowed with a finite dimensional commutative subalgebra $\D$. Note that such a situation can arise if one considers operators between different spaces, say $H_1, H_2$, an embeds them in $B(H_1\oplus H_2)$ as it was made for matrices in (\ref{je.chancelle.un.tt.petit.peu.quand.je.te.rencontre}). We define a microstate free entropy for $N$-tuples $(a_1,\ldots, a_N)$ of elements of $\A$  which have kernel and range projections in $\D$: it is the asymptotic logarithm of the volume of $N$-tuples of rectangular matrices whose joint distribution (under the state defined by the trace) is closed to  the joint distribution of $(a_1,\ldots, a_N)$ in $(\A,\vfi)$. 

This free entropy is subadditive, and additive only on families which are free with amalgamation over $\D$. This is one of the properties that has made us consider this free entropy possibly relevant to study operators between different Hilbert spaces with the tools of free probability. 

Another optimization problem has given rise to an interesting analogy.  In the previous paper \cite{fbg.AOP.rect}, for each $\la\in (0,1)$, we defined a free convolution $\mu\arc\nu$ of symmetric \pro measures as the   distribution in $(q\A q, \ff{\vfi(q)}\vfi)$ of $a+b$, where  $a, b$ are free with amalgamation over $\D$, have kernel projection $\leq p=1-q\in \D$ and range projection $\leq q$ \st $\f{\vfi(q)}{\vfi(p)}=\la$, and have symmetrized distributions $\mu,\nu$ in $(q\A q, \ff{\vfi(q)}\vfi)$. We established in \cite{fbg05.inf.div} a correspondence (like the Bercovici-Pata bijection) between  $\arc$-infinitely divisible distributions and $*$-infinitely divisible distributions. In this correspondence, the analogue of Gaussian distributions are symmetrizations of Marchenko-Pastur distributions. In this paper, we prove that among the set of elements $a$ with kernel projection $\leq p$ and range projection $\leq q$ and \st $\vfi(aa^*)\leq \vfi(q)$, the elements which maximize free entropy are the elements $a$ \st the distribution of $aa^*$ in $(q\A q, \ff{\vfi(q)}\vfi)$ is a Marchenko-Pastur distribution. 

We also construct  a free Fisher's information with the conjugate variables approach for elements which have kernel and range projections in $\D$. We have a Cram\'er-Rao inequality, where Marchenko-Pastur distributions appear again as the distributions which realize equality, and  a superadditivity result where freeness with amalgamation over $\D$ is equivalent to additivity (when quantities are finite).

The main relevance, according to the author, of this problems of optimization, is the {\it legitimization} of this notions. Indeed, the analoguous problems for the classical entropy and information in one hand, and for the  entropy and the information defined by Voiculescu one the other hand, have been solved (see \cite{HP99}, \cite{hiai}, \cite{NSS99}, \cite{NSS99.2},  \cite{spei99}), 	and the solutions where actually the analogues of the solutions given here. This supports the idea that the notions proposed here are the right ones to apply the tools and the ideas of free probability theory to the study of operators between different Hilbert spaces. Moreover, the solutions of optimization problems for entropy and information under certain constraints are, in a sens, the {\it generic} objects which realize this constraints. 

In section 1 and 2, we define the objects we are going to use and we recall definitions and basic properties of operator valued cumulants.

In section 3, we prove that under certain hypothesis, freeness with respect to the state $\vfi$ implies freeness with amalgamation over the finite dimensional commutative algebra $\D$. As an application, we prove a result about polar decomposition of $R$-diagonal elements with non trivial kernel: the partial isometry and the positive part are  free with amalgamation over the algebra generated by the kernel projection. 

In section 4, we prove asymptotic freeness with amalgamation over $\D$ of rectangular independent random matrices (as a consequence of results of the previous section). This result is used section 5, where we define our microstates free entropy and solve the optimization problems we talked about above, using some change of variable formulae we establish in the same section. Similarly, in section 6, we construct our free Fisher's information with the conjugate variables approach and solve optimization problems.

{\bf Aknowledgements.} We would like to thank Philippe Biane, Dan Voiculescu, and Piotr \'Sniady for useful discussions, as well as Thierry Cabanal-Duvillard, who organized the workshop ``Journ\'ee Probabilit\'es Libres'' at MAP5 in June 2004, where the author had the opportunity to have some of these discussions.

\section{Definitions}\label{amis+frere.26.07.05}
In this section, we will define the spaces and the notions. For all $d$ integer, we denote by $[d]$ the set $\{1,\ldots, d\}$.

Consider a  tracial $*$-noncommutative \pro space $(\mc{A},\vfi)$  endowed with a family $(p_1,\ldots,p_d)$ of self-adjoint non zero projectors (i.e. $\forall i, p_i^2=p_i$) which are pairwise orthogonal (i.e. $\forall i\neq j, p_ip_j=0$), and such that $p_1+\cdots +p_d=1$. Any element $x$ of $\mc{A}$ can then be represented $$x=\begin{bmatrix}x_{11}& \cdots & x_{1d}\\  \vdots & & \vdots \\ x_{d1}& \cdots & x_{dd}\end{bmatrix},$$ where $\forall i,j, x_{ij}=p_ixp_j$. This notation is compatible with the product and the involution. 

Let us define, for all $i,j\in [d]$, $\mc{A}_{i,j}=p_i\mc{A}p_j$ (the comma between $i$ and $j$ will often be omitted). We call {\it simple elements} the non zero elements of the union of the $\mc{A}_{ij}$'s ($i,j\in [d]$).  We define $\vfi_i:=\ff{\rho_i}\vfi_{|\mc{A}_ {ii}}$, with $\rho_i:=\vfi(p_i)$. Note that, since $\vfi$ is a trace, every $\vfi_i$ is a trace, but for  $i, j\in[d]$, $a\in\mc{A}_{ij}$, $b\in \mc{A}_{ji}$, one has \begin{equation}\label{mai68,agelyrique?}\rho_i\vfi_i(ab)=\rho_j\vfi_j(ba).\end{equation} 

Note also that the linear span $\mc{D}$ of $\{p_1,\ldots,p_d\}$ is a $*$-algebra, which will be identified to the set of $d\tii d$ diagonal complex matrices by $$\sum_{i=1}^d\la_i p_i\simeq \diag (\la_1,\ldots,\la_d).$$ The application $\ED$, which maps $x\in \mc{A}$ to $\diag(\vfi_1(x_{11}),\ldots , \vfi_d(x_{dd}))$, is then a conditional expectation from $\mc{A}$ to $\D$: $$\forall (d,a,d')\in \mc{D}\times \mc{A}\times \mc{D}, \ED(dad')=d\ED(a)d'.$$

 A family $(\mc{A}_i)_{i\in I}$ of subalgebras of $\mc{A}$ which all contain $\mc{D}$ is said to be {\it free with amalgamation over $\mc{D}$} 
if for all $n$, $i_1\neq\cdots\neq i_n\in I$, for all $x(1)\in \mc{A}_{i_1}\cap \ker \ED, \ldots, x(n)\in \mc{A}_{i_n}\cap \ker \ED$, one has \begin{equation}\label{atomheartmother}\ED(x(1))\cdots x(n))=0.\end{equation} A family $(\chi_i)_{i\in I}$ of subsets of $\A$ is said to be free with amalgamation over $\D$ if there exists free with amalgamation over $\D$ subalgebras $(\mc{A}_i)_{i\in I}$ (which all contain $\mc{D}$) \st for all $i$, $\chi_i\subset \A_i$.

The {\it $\mc{D}$-distribution} of a family $(a_i)_{i\in I}$ of elements of $\mc{A}$ is the application which maps a word $X_{i_0}^{\epsilon_0}d_1X_{i_1}^{\epsilon_1}\cdots X_{i_n}^{\epsilon_n}d_n$ in $X_i,X_i^*$ ($i\in I$) and the elements of $\D$ to  $\ED (a_{i_0}^{\epsilon_0}d_1a_{i_1}^{\epsilon_1}\cdots a_{i_n}^{\epsilon_n}d_n)$. 

It is easy to see that the $\mc{D}$-distribution of a free with amalgamation over $\D$  family depends only on the individual $\D$-distributions.

Consider a sequence $(\A_n,\vfi_n)$ of tracial $*$-noncommutative \pro spaces \st for all $n$, $\D$ can be identified with a $*$-subalgebra of $\A_n$ (the identification is not supposed to preserve the state). 
The {\it convergence in $\D$-distribution} of a sequence $(a_i(n))_{i\in I}$ of families of elements of the $\A_n$'s to a family $(a_i)_{i\in I}$ of $\A$ is the pointwise convergence of the sequence of $\D$-distributions. In this case, if $I=\cup_{s\in S}I_s$ is  a partition of $I$, then the family of subsets $(\{a_i(n)\ste i\in I_s\})_{s\in S}$ is said to be {\it asymptotically \faose} if the family of subsets $(\{a_i\ste i\in I_s\})_{s\in S}$ is \faop

It is easy to see that the $\mc{D}$-distribution of a free with amalgamation over $\D$  family depends only on the individual $\D$-distributions.

\section{Cumulants}
The theory of cumulants in a $\D$-\pro space (i.e. in an algebra endowed with a conditional expectation on a subalgebra $\D$) has been developed in \cite{spei98}. In this section, we will begin by giving the main lines of this theory, and then we will investigate the special case of the situation we presented in the previous section.

\subsection{General theory of cumulants in a $\D$-\pro space}\label{friendsperdu}
In this section, we consider an algebra $\mc{A}$, a subalgebra $\D$ of $\mc{A}$, and a conditional expectation  $\ED$ form $\mc{A}$ to $\D$. 
\\
Let us begin with algebraic definitions. A $\D$-bimodule is a vector space $M$ over $\C$ on which the algebra $\D$ acts on the right and on the left. 
The tensor product $M\otimes_\D N$ of two  $\D$-bimodules $M,N$ is their tensor product as $\C$-vector spaces, where for all $(m,d,n)\in M\stt \D\stt N$, $(m.d)\otimes n$ and $m\otimes (d.n)$ are identified. $ M\otimes_\D N$ is endowed with a structure of  $\D$-bimodule by $d_1.(m\otimes n).d_2=(d_1.m)\otimes (n.d_2)$. This allows us to define, for $n$ positive integer, $\A^{\otimes_\D n}=\underbrace{\A\otimes_\D\cdots \otimes_\D\A}_{\textrm{$n$ times}}$.
\\
Consider a sequence $(f_n)_{n\geq 1}$ of maps, each $f_n$ being a $\D$-bimodule morphism between $\A^{\otimes_\D n}$ and $\D$. For $n$ positive integer and $\pi\in \NC(n)$ (noncrossing partition of $[n]$), we define the   $\D$-bimodule morphism $f_\pi$ between $\A^{\otimes_\D n}$ and $\D$  in the following way: if $\pi=1_n$ is the one-block partition, $f_\pi=f_n$. In the other case, a block $V$ of $\pi$ is an interval $[k,l]$. If $k=1$ (resp. $l=n$), then $f_\pi(a_1\otimes\cdots\otimes a_n)=f_{l-k+1}(a_1\otimes\cdots \otimes a_l)f_{\pi\backslash \{V\}}(a_{l+1}\otimes\cdots\otimes a_n)$ (resp. $f_{\pi\backslash \{V\}}(a_{1}\otimes\cdots\otimes a_{k-1})f_{l-k+1}(a_k\otimes\cdots \otimes a_n)$). In the other case, one has $1<k\leq l<n$. Then $ f_\pi(a_1\otimes\cdots\otimes a_n)$ is defined to be $f_{\pi\backslash \{V\}}(a_{1}\otimes\cdots\otimes a_{k-1}f_{l-k+1}(a_k\otimes\cdots \otimes a_l)\otimes a_{l+1}\otimes\cdots\otimes a_n)$ or  $f_{\pi\backslash \{V\}}(a_{1}\otimes\cdots\otimes a_{k-1}\otimes f_{l-k+1}(a_k\otimes\cdots \otimes a_l)a_{l+1}\otimes\cdots\otimes a_n)$, both are the same by definition of $\otimes_\D$. 
\\
For example, if $\pi=\{\{1,6,8\},\{2,5\},\{3,4\},\{7\},\{9\}\}$, then $$f_\pi(a_1\otimes\cdots\otimes a_9)=f_3\lf(a_1f_2\lf(a_2f_2(a_3\otimes a_4)\otimes a_5\ri)\otimes a_6f_1(a_7)\otimes a_8\ri)f_1(a_9).$$
\\
Let us define, for all $n\geq 1$, the $\D$-bimodule morphism $\EDn$ between $\A^{\otimes_\D n}$ and $\D$ which maps $a_1\otimes\cdots\otimes a_n$ to $\ED(a_1\cdots a_n)$. Then one can define the sequence $(c_n)_{n\geq 1}$ of maps, each $c_n$ being a $\D$-bimodule morphism between $\A^{\otimes_\D n}$ and $\D$,
by one of the following equivalent formulae:
\begin{eqnarray}\ds\forall n, \forall \pi\in \NC(n),\quad\quad  \EDpi&=&\sum_{\substack{\sigma\leq \pi}}c_\sigma,\label{def.cum.1.1.05.1}\\
\ds\forall n,\quad\quad \EDn &=& \sum_{\substack{\sigma\in \NC(n)}}c_\sigma\label{def.cum.1.1.05.2}\\
\ds\forall n, \forall \pi\in \NC(n), \quad\quad c_\pi&=&\sum_{\substack{\sigma\leq \pi}}\mu(\sigma,\pi)\EDsigma,\label{def.cum.1.1.05.3}\\
\ds\forall n, \quad\quad c_n& =& \sum_{\substack{\sigma\in \NC(n)}}\mu(\sigma,1_n)\EDsigma,\label{def.cum.1.1.05.4}
\end{eqnarray}
where $\mu$ is the M\"obius function  (\cite{R64}, \cite{spei99}) of the lattice $\NC(n)$ endowed with the reffinment order. 
\\
The following result is a consequence of Proposition  3.3.3 of \cite{spei98}, used with the formula of cumulants with products as entries (Theorem 2 of \cite{SS01}), which can be generalized to $\D$-\pro spaces. 
\begin{Th}\label{article.tombant.en.ruine}  A family $(\chi_i)_{i\in I}$ of subsets of $\A$ is \fao \ssi for all $n\geq 2$, for all non constant $i\in I^n$, for all $a_1\in \chi_{i_1}$,...,   $a_n\in \chi_{i_n}$, one has $c_n(a_1\otimes\cdots\otimes a_n)=0$.
\end{Th}
Note that this \teo is a little improvement of  Theorem 1 of \cite{SS01}.

\subsection{The special case where $\D=\Span (p_1,\ldots,p_d)$}\label{13.7.4.1}
For the rest of the text, we consider again, without introducing them, the same objects  as in section \ref{amis+frere.26.07.05}. By linearity of the cumulant functions, we will work only with  simple elements (i.e. non zero elements of the union of the $\mc{A}_{ij}$'s, $1\leq i,j\leq d$).

(a) First, for all $i,j,k,l\in [d]$ \st $j\neq k$, one has $\A_{ij}\otimes_\D \A_{kl}=\{0\}$ (because $p_jp_k=0$). 
So we will only have to compute the cumulant functions on subspaces of the type $\A_{i_0i_1}\otimes_\D \A{i_1i_2}\otimes_\D\cdots \otimes_\D \A{i_{n-1}i_n}$, with $i_0,i_1,\ldots,i_n\in [d]$.

(b) Moreover, on such a subspace, $c_n$ takes values in $\A_{i_0i_n}$, because it is a $\D$-bimodule morphism. So, if $i_0\neq i_n$, since $\D\cap\A_{i_0i_n}=\{0\}$, $c_n$ is null on  $\A_{i_0i_1}\otimes_\D \A_{i_1i_2}\otimes_\D\cdots \otimes_\D \A_{i_{n-1}i_n}$.
So it is easily proved by induction that  for $\pi\in \NC(n)$,  for all $i_0,i_1$,..., $i_n\in [d]$, $c_\pi$ is null on $\A_{i_0i_1}\otimes_\D \A{i_1i_2}\otimes_\D\cdots \otimes_\D \A{i_{n-1}i_n}$ whenever a block $\{k_1<$ ... $ <k_m\}$ of $\pi$ is \st $i_{k_1-1}\neq i_{k_m}$.   

(c) Hence the function $c_\pi$ factorizes on the complex vector space $\A_{i_0i_1}\otimes_\D \A_{i_1i_2}\otimes_\D\cdots \otimes_\D \A_{i_{n-1}i_n}$ in the following way: for $(a_1,$ ... $,a_n )\in \A_{i_0i_1}\times\cdots\times\A_{i_{n-1}i_n}$,  
\begin{equation}\label{volanniv}\ds c_\pi(a_1\otimes \cdots \otimes a_n)=\lf(\prod_{\substack{V\in \pi\\ V=\{k_1<\cdots <k_m\}}}c^{(i_{k_m})}_m(a_{k_1}\otimes\cdots\otimes a_{k_m})\ri).p_{i_n}, \end{equation}
where for all $m$, $c^{(1)}_m$,..., $c^{(d)}_m$ are the linear forms  on the complex   vector space $\A^{\otimes_\D m}$ defined by $$\forall x\in\A^{\otimes_\D m}, c_m(x)=\begin{bmatrix}c_m^{(1)}(x)& & \\ & \ddots & \\ & & c_m^{(d)}(x)\end{bmatrix}.$$
Formula (\ref{volanniv}) can be written in the following way:  for $(a_1,$ ... $,a_n )\in \A_{i_0i_1}\times\cdots\times\A_{i_{n-1}i_n}$, 
\begin{equation}\label{volanniv.2.3.05.1}\ds c_\pi(a_1\otimes \cdots \otimes a_n)=\prod_{\substack{V\in \pi\\ V=\{k_1<\cdots <k_m\}}}\eta_{i_n, i_{k_m}}\circ c_m(a_{k_1}\otimes\cdots\otimes a_{k_m}), \end{equation}
where for all $i,j\in [d]$, $\eta_{i,j}$ is the involution of $\D$ which permutes the $i$-th and the $j$-th columns in the representation of elements of $\D$ as $d\tii d$ complex matrices.

\begin{rmq}\label{14.7.4.1} In (b), (c), we only used the fact that for all $n$, $c_n$ is a $\D$-bimodule morphism, so everything stays true if one replaces $c_n$ by  $\EDn$ and $c_\pi$ by $\EDpi$.\end{rmq}

 (d) Now it remains only to investigate the relation between the functions $c_n^{(1)}$,...,  $c_n^{(d)}$.  We will prove, by induction on $n$, a formula analogous to (\ref{mai68,agelyrique?}). Consider $(a_1,$ ... $,a_n )\in \A_{i_0i_1}\times\cdots\times\A_{i_{n-1}i_n}$, with $i_0=i_n$. Then one has \begin{equation}\label{ecologie&feminisme}\rho_{i_0} c_n^{(i_0)}(a_1\otimes \cdots\otimes a_n)=\rho_{i_1}c_n^{(i_1)}(a_2\otimes \cdots\otimes a_n\otimes a_1).\end{equation}
For $n=1$, it is clear. Now suppose the result proved to the ranks $1,\ldots,n-1$, and consider $(a_1,$ ... $,a_n )\in \A_{i_0i_1}\times\cdots\times\A_{i_{n-1}i_n}$, with $i_0=i_n$. One has, by formulae (\ref{def.cum.1.1.05.2}),(\ref{volanniv}),\begin{equation*}c_n^{(i_0)}(a_1\otimes \cdots\otimes a_n)=\underbrace{\vfi_{i_0}(a_1\cdots a_n)}_{X}-\underbrace{\sum_{\substack{\pi\in \NC(n)\\ \pi<1_n}}\prod_{\substack{V\in \pi\\ V=\{k_1<\cdots <k_m\}}}c^{(i_{k_m})}_m(a_{k_1}\otimes\cdots\otimes a_{k_m})}_{Y},\end{equation*} and  \begin{equation*}c_n^{(i_1)}(a_2\otimes \cdots\otimes a_n\otimes a_1)=\underbrace{\vfi_{i_1}(a_2\cdots a_na_1)}_{X'}-\underbrace{\sum_{\substack{\pi\in \NC(n)\\ \pi<1_n}}\prod_{\substack{V\in \pi\\ V=\{k_1<\cdots <k_m\}}}c^{(i_{\sigma(k_m)})}_m(a_{\sigma(k_1)}\otimes\cdots\otimes a_{\sigma(k_m)})}_{Y'},\end{equation*} where $\sigma$ is the cycle $(12\cdots n)$ of $[n]$.
\\
Since $\rho_{i_0} X=\rho_{i_1}X'$ (by formula (\ref{mai68,agelyrique?})), it suffices to prove that   \begin{equation*}\rho_{i_0}Y=\rho_{i_1}Y'.\end{equation*} To do that, it suffices to propose a bijective correspondence $\pi\mapsto \tilde{\pi}$ form $\NC(n)-\{1_n\}$ to $\NC(n)-\{1_n\}$ such that for all $\pi\in \NC(n)-\{1_n\}$, $$\rho_{i_0}\prod_{\substack{V\in \pi\\ V=\{k_1<\cdots <k_m\}}}c^{(i_{k_m})}_m(a_{k_1}\otimes\cdots\otimes a_{k_m})=\rho_{i_1} \prod_{\substack{V\in \tilde{\pi}\\ V=\{k_1<\cdots <k_m\}}}c^{(i_{\sigma(k_m)})}_m(a_{\sigma(k_1)}\otimes\cdots\otimes a_{\sigma(k_m)}).$$ By induction hypothesis, the correspondence which maps $\pi\in\NC(n)-\{1_n\}$ to $\tilde{\pi}$ defined by $$k\stackrel{\tilde{\pi}}{\sim} l\Leftrightarrow \sigma(k)\stackrel{\pi}{\sim}\sigma( l)$$ is convenient.

The following \teo has been proved in the section called {\it Rectangular Gaussian distribution and Marchenko-Pastur distribution} of \cite{fbg05.inf.div}.
\begin{Th}\label{regard.qui.perce.et.qui.met.quelque.chose.en.feu.dans.moi} For $k,l\in [d]$ \st $\rho_k\leq \rho_l$,  $b\in \A_{k,l}$ satisfies, for all positive integer $n$, $$c_{2n}^{(k)}(b\otimes b^*\otimes\cdots \otimes b^*)=\f{\rho_l}{\rho_k}\delta_{n,1} $$\ssi the moments of $bb^*$ in $(\A_{kk},\vfi_k)$ is the Marchenko-Pastur distribution with parameter $\f{\rho_l}{\rho_k}$ (defined p. 101 of \cite{hiai}).
\end{Th}

\section{Freness with respect to $\vfi$ versus freeness with amalgamation over $\D$} 
\subsection{$\D$-central limit theorems} On sets of matrices, $||.||$ will denote the operator norm associated to the canonical hermitian norms.
A self-adjoint element $X$ of $\A$ is said to be {\it $\D$-semicircular with covariance $\vfi$} if it satisfies \begin{itemize}\item[(i)] $c_1(X)=0$
\item[(ii)] $\forall d\in\D, \;\; c_2(Xd\otimes X)=\vfi(d)$,
\item[(iii)] $\forall k\geq 3, \forall d_1,\ldots, d_{k}\in\D,\;\; c_k(Xd_1\otimes\cdots\otimes Xd_{k})=0$ 
\end{itemize}
Note that it determines the $\D$-distribution of $X$.

\begin{Th}[$\D$-central limit theorem] Consider a family $(X_i)_{i\geq 1}$ of self-adjoint elements of $\A$ which satisfy\begin{itemize}\item[(a)] $X_1,X_2,$... are \faose,\item[(b)] $\forall i,\forall d\in \D, \ED(X_i)=0,\ED(X_idX_i)=\vfi(D)$, \item[(c)] $\forall k,\forall d_1$,..., $d_{k}\in \D, \sup_{i\geq 1}||\ED(X_id_1\cdots X_id_k)||<\infty.$   \end{itemize} Then $Y_n:=\ff{\sqrt{n}}\sum_{i=1}^n X_i$ converges in $\D$-distribution to a $\D$-semicircular element with covariance $\vfi$.
\end{Th}

This theorem is very closed to many well-known results of free probability theory (e.g. Theorem 4.2.4 of \cite{spei98}). 

 We prove now a kind of multidimentional $\D$-central limit Theorem, analoguous to Theorem 2.1 of \cite{voic91}:
\begin{Th}\label{theo.alg.noel.04} Consider a family $(T_j)_{j\in\n}$ of self-adjoint elements of $\A$, a $*$-subalgebra $\B$ of $\A$ containing $\D$, \st \\
(H1) $\forall m,\forall B_1$,..., $B_{m}\in \B, \sup_{i_1,\ldots,i_m\in\n}||\ED(T_{i_1}B_1T_{i_2}\cdots T_{i_m}B_m)||<\infty,$
\\
(H2) for $m\geq 1$, for $B_0$,..., $B_m\in\B$, for $\alpha : [m]\to \n$,  
one has  \begin{itemize}  
\item[(a)] $\ED(B_0T_{\alpha(1)}B_1\cdots T_{\alpha(m)}B_m)=0$ if an element of $\n$ has exactly one antecedent by $\alpha$,
\item[(b)]  $\ED(B_0T_{\alpha(1)}B_1\cdots T_{\alpha(m)}B_m)=\vfi(B_r) \ED(B_0T_{\alpha(1)}B_1\cdots T_{\alpha(r-1)}B_{r-1}B_{r+1}\cdots T_{\alpha(m)}B_m)$ if no element of $\n$ has strictly more than two antecedents by $\alpha$ and $\alpha(r)=\alpha(r+1)$, with $1\leq r<m$, 
\item[(c)] $\ED(B_0T_{\alpha(1)}B_1\cdots T_{\alpha(m)}B_m)=0$  if no element of $\n$ has strictly more than two antecedents by $\alpha$ and for all  $1\leq r<m$, $\alpha(r)\neq\alpha(r+1)$.
\end{itemize}
Consider $\beta :\n^2\to\n$ injective, and define $X_{m,n}=\ff{\sqrt{n}}\sum_{j=1}^nT_{\beta(m,j)}$. Then for all $m$, $X_{m,n}$ converge in distribution, when $n\to\infty$, to a  $\D$-semicircular element with covariance $\vfi$, and the family of subsets $(B, (\{X_{m,n}\})_{m\in\n})$ is asymptotically \fao as $n\to\infty$. \\
Moreover, if 
\begin{itemize}\item[(a')] $(B, (\{T_{j}\})_{j\in\n})$ is \faose, 
\item[(b')] $\forall j,\forall d\in \D, \ED(T_j)=0,\ED(T_jdT_j)=\vfi(d)$,
\item[(c')] $\forall m,\forall d_1$,..., $d_{m}\in \D, \sup_{j\in\n}||\ED(T_{j}d_1\cdots T_{j}d_m)||<\infty,$
\end{itemize} then (H1) and (H2) are satisfied. 
\end{Th}

\begin{pr}We shall proced as in the proof of theorem 2.1 of \cite{voic91}. First we prove that [(a'),(b'),(c')] implies [(H1),(H2)]. Then  we prove that it suffices to prove the result replacing  [(H1),(H2)] by  [(a'),(b'),(c')], and at last we prove the result in this particular case. For $x\in \A$, we define $\overset{\circ}{x}:=x-\ED(x)$.

{\it Step I. } Suppose that the $T_j$'s and $\D$ satisfy [(a'),(b'),(c')].

The proof of the fact that (a') and (c') together implie (H1) is along the same lines as the proof of $1\,^{\circ}$ of the Step I of the proof of Theorem 2.1 of \cite{voic91}, so we leave it to the reader.

Consider  $m\geq 1$,  $B_0$,..., $B_m\in\B$, $\alpha : [m]\to \n$.

(H2).(a) follows from (a'), (b'), and the following easy result: \begin{equation}\label{referenceexcellente}
\forall a,b,c\in\A, [\{a\},\{b,c\}\textrm{\faose}]\Rightarrow \ED(bac)=\ED(b\ED(a)c).\end{equation}

Suppose no element of $\n$ has strictly more than two antecedents by $\alpha$ and $\alpha(r)=\alpha(r+1)$, with $1\leq r<m$. 
\\
Let us prove   $\ED(B_0T_{\alpha(1)}B_1\cdots T_{\alpha(m)}B_m)=\vfi(B_r) \ED(\underbrace{B_0T_{\alpha(1)}B_1\cdots T_{\alpha(r-1)}B_{r-1}}_{:=A}\underbrace{B_{r+1}\cdots T_{\alpha(m)}B_m}_{:=B})$. 
\\
Suppose first that $B_r\in\D$. Then $\{Y_{\alpha(r)}B_rY_{\alpha(r+1)}\}$, $\{A,B\}$ are \faose, so, by (\ref{referenceexcellente}),  
$$\ED(B_0T_{\alpha(1)}B_1\cdots T_{\alpha(m)}B_m)=\ED(A\ED(Y_{\alpha(r)}B_rY_{\alpha(r+1)})B).$$
 But by (b'), $\ED(Y_{\alpha(r)}B_rY_{\alpha(r+1)})=\vfi(B_r)$, which allows us to conclude.
\\ So, by linearity, we can now suppose that $\ED(B_r)=0$. In this case, $\vfi(B_r)=0$, so it suffices to prove that  $\ED(B_0T_{\alpha(1)}B_1\cdots T_{\alpha(m)}B_m)=0$. It follows from (a') and (\ref{atomheartmother}), applied to all terms of the right hand side of: \begin{eqnarray*}\ED(B_0T_{\alpha(1)}B_1\cdots T_{\alpha(m)}B_m)&=&\ED(\overset{\circ}{A}Y_{\alpha(r)}B_rY_{\alpha(r+1)}\overset{\circ}{B})+\ED(A)\ED(Y_{\alpha(r)}B_rY_{\alpha(r+1)}\overset{\circ}{B})\\ &&+\ED(\overset{\circ}{A}Y_{\alpha(r)}B_rY_{\alpha(r+1)})\ED(B)+\ED(A)\ED(Y_{\alpha(r)}B_rY_{\alpha(r+1)})\ED(B).
\end{eqnarray*}

Suppose that no element of $\n$ has strictly more than two antecedents by $\alpha$ and that for all $1\leq r<m$, $\alpha(r)\neq\alpha(r+1)$. By linearity,  $\ED(B_0T_{\alpha(1)}B_1\cdots T_{\alpha(m)}B_m)$ is equal to 
$$\sum_{P\subset\{0,\ldots,m\}}\ED\lf((1_P(0)\overset{\circ}{B_0}+1_{P^c}(0)\ED(B_0))T_{\alpha(0)}\cdots T_{\alpha(m)}(1_P(m)\overset{\circ}{B_0}+1_{P^c}(m)\ED(B_0))\ri),$$ where for $P\subset\{0,\ldots,m\}$, $1_P$ (resp. $1_{P^c}$) denotes the characteristic function of $P$ (resp. of its complementary).  It follows from (a') and (\ref{atomheartmother}), applied to all terms of the sum, that $\ED(B_0T_{\alpha(1)}B_1\cdots T_{\alpha(m)}B_m)=0$.

{\it Step II. } After having eventually extended $\A$, consider  a \fao family $(x_m)_{m\geq 1}$  of  $\D$-semicircular elements of $\A$ with covariance $\vfi$, which is also \fao with $\B$. Let us show that in order to prove that for all $r\geq 1$, $B_0$,..., $B_r\in \B$, $m : [r]\to\n$, 
$$\ED(B_0X_{m(1),n}B_1\cdots B_{r-1}X_{m(r),n}B_r)\underset{n\to\infty}{\longrightarrow}\ED(B_0 x_{m(1)}B_1\cdots B_{r-1}x_{m(r)}B_r),$$ it suffices to prove it in the particular case where [(a'),(b'),(c')] are satisfied.

So consider $r\geq 1$, $B_0$,..., $B_r\in \B$, and $m : [r]\to\n$. Define, for $n\geq 1$, the set $P_n=m([r])\tii [n]$, and define, for $I=(p_1$,..., $p_r)\in P_n^r$, $$\Pi_I=B_0T_{\beta(p_1)}B_1\cdots B_{r-1}T_{\beta(p_r)}B_r.$$ Then by linearity, there exists a family $(C_I)_I$ of elements of $\{0,1\}$, indexed by $I\in P_n^r$, \st we have:
$$\ED(B_0X_{m(1),n}B_1\cdots B_{r-1}X_{m(r),n}B_r)=\ff{n^{\f{r}{2}}}\sum_{I\in P_n^r}C_I\ED(\Pi_I).$$ 
\\
By (H2).(a), if $\ED(\Pi_I)\neq 0$, then no element of $P_n$ appears exactly once in $I$. Let $R_{n,r}$ be the set of elements $I$ of $P_n^r$ \st no element of $P_n$ appears exactly once in $I$ and an element of $P_n$ appears at least three times in $I$. Its cardinality is less than $|P_n|\tii |P_n|^{(r-3)/2}r!=o(n^{\f{r}{2}})$, so, since by (H1) there exists $M>0$ \st for all $n$, $I\in P_r^n$, $|| \ED(\Pi_I)||\leq M$, one has 
$$\ff{n^{\f{r}{2}}}\sum_{I\in R_{n,r}}||C_I\ED(\Pi_I)||\underset{n\to\infty}{\longrightarrow}0.$$
\\
So \begin{equation}\label{nahum.family}\ds\lim_{n\to\infty}\ED(B_0X_{m(1),n}B_1\cdots B_{r-1}X_{m(r),n}B_r)\end{equation} exists \ssi $$\ds\lim_{n\to\infty}\ff{n^{ \f{r}{2} }}\sum_{\substack{I\in P_n^r,\\ \textrm{each element of $P_r$} \\ \textrm{appears exactly}\\ \textrm{$0$ or $2$ times in $I$}}}C_I\ED(\Pi_I),$$ exists, and in this case, the limits are the same. 
\\
But the computation of $\ED(\Pi_I)$, for elements $I$ of $P_n^r$ such as those considered in the previous sum, is completely determined by (H2). So the limit (\ref{nahum.family}) will be the same (and exist in the same time) if one replaces the $T_j$'s by another family which satisfies [(H1),(H2)]. In particular, by Step I, one can suppose that [(a'),(b'),(c')] are satisfied.

{\it Step III. } Suppose now that  [(a'),(b'),(c')] are satisfied. The previous theorem allows us to claim that for all $m$, $X_{m,n}$  converges in $\D$-distribution, as $n\to\infty$, to a $\D$-semicircular with covariance $\vfi$. Moreover, for all $n$,  the family of subsets $(B, (\{X_{m,n}\})_{m\in\n})$ is \faose, so the theorem is proved. 
\end{pr}


The main theorem of this section is the following one. Recall that a family $(\A_{i})_{i\in I}$ of subalgebras of $\A$ is said to be {\it free} if for all $n$, $i_1\neq\cdots\neq i_n\in I$, for all $x(1)\in \mc{A}_{i_1}\cap \ker \vfi, \ldots, x(n)\in \mc{A}_{i_n}\cap \ker \vfi$, one has \begin{equation}\label{atomheartmother.2.2.05}\vfi(x(1)\cdots x(n))=0.\end{equation} A family $(\chi_i)_{i\in I}$ of subsets of $\A$ is said to be {\it free} if there exists free subalgebras $(\mc{A}_i)_{i\in I}$ \st for all $i$, $\chi_i\subset \A_i$. In order to avoid confusion between freeness and freeness with amalgamation over $\D$, freeness will be called {\it $\vfi$-freeness}. 
We use the notion of {\it $\vfi$-distribution} of a family $(a_i)_{i\in I}$ of elements of $\A$: it is the application which maps a word   $X_{i_1}^{\epsilon_1}\cdots X_{i_n}^{\epsilon_n}$ in $X_i,X_i^*$ ($i\in I$) to  $\vfi (a_{i_1}^{\epsilon_1}\cdots a_{i_n}^{\epsilon_n})$. It is easy to see that the $\vfi$-distribution of a $\vfi$-free  family depends only on the individual $\vfi$-distributions, and that the $\D$-distribution of a family which contains $\D$ is determined by its $\vfi$-distribution. At last, recall that $\vfi$-semicircular elements are elements whose moments are given by the moments of the semicircle distribution with center $0$ and radius $2$.

\begin{Th}\label{3.1.05.1} Consider, in $\A$, a family $(y(s))_{s\in \n}$ of $\vfi$-semicircular elements, 
and a subalgebra $\B$ of $\A$ which contains $\D$ \st the family $(\B, (\{y(s)\})_{s\in\n})$ is $\vfi$-free. Then    the family $(\B, (\{y(s)\})_{s\in\n})$ is also free with amalgamation over $\D$, and the $\D$-distribution of $y(s)$'s is the $\D$-semicircular  distribution   with covariance $\vfi$.
\end{Th}

\begin{pr} Consider $\beta:\n\tii\n\to\n$ injective. By stability of $\vfi$-semicircular distribution under free convolution, it is clear that for all $n\geq 1$, the family $$(\B,\ff{\sqrt{n}}\sum_{j=1}^n  y(\beta(m,j)))_{m\geq 0})$$ has the same $\vfi$-distribution (and hence $\D$-distribution, because contains $\D$) as $(\B, (\{y(s)\})_{s\in\n})$. So it suffices to prove that $(\B, (\{y(s)\})_{s\in\n})$ satisfies (H1) and (H2).

(H1) is due to the fact that for $m\geq 1$ and $ b_1,\ldots, b_m\in\B$ fixed,  \begin{eqnarray*}\ds \forall s_1,\ldots, s_m\in\n,
||\ED(  y(s_1)b_1\cdots y(s_m)b_m)||&=&
|| \sum_{k=1}^d \vfi_k(p_ky(s_1)b_1\cdots y(s_m)b_mp_k).p_k||\\
\ds & =&\max_{1\leq k\leq d}\ff{\rho_k}\lf|\vfi(p_ky(s_1)b_1\cdots y(s_m)b_mp_k)\ri|,\end{eqnarray*}
which only depends on the partition $\pi$ of $[m]$ which links two elements $i,j$ \ssi $s_i=s_j$.

To prove (H2)(a),(b),(c), since for all $x\in \A$, $$\ED(x)=\sum_{k=1}^d\ff{\rho_i}\vfi(p_ixp_i).p_i$$ and the algebra $\B$ contains all $p_i$'s,    it suffices to prove it with $\ED$ replaced by $\vfi$. Then it follows from the last assertion of Theorem   2.1 of \cite{voic91}. 
\end{pr}

\subsubsection*{Remarks about the previous theorem}
\rm (1) Let $\mc{C}$ be a subalgebra of $\A$ which is $\vfi$-free with $\D$. It is easy to see that for all $x\in \mc{C}$, $\ED(x)=\vfi(x).1$, hence for all $n\geq 1$, $a_1,\ldots, a_n\in\mc{C}$, $$c_n(a_1\otimes\cdots\otimes a_n)=\mf{K}_n(a_1,\ldots,a_n).1, $$ where $\mf{K}_n$ is the $n$-th $\vfi$-cumulant function. 
So a $\vfi$-free family of subalgebras of $\mc{C}$ is also $\ED$-free: $\vfi$-freeness implies vanishing of mixed $\vfi$-cumulants, which implies $\ED$-freeness, by theorem \ref{article.tombant.en.ruine}. It is not enough to prove our result, because the algebra $\mc{C}$ cannot in the same time be $\vfi$-free with $\D$ and  contain $\D$, hence cannot contain $\B$. 
\\
(2) This theorem recalls Theorem 3.5 of \cite{NSS02}. But to prove our result using this theorem, it would be necessary to compute $\B$-cumulant functions. 

\subsection{Polar decomposition $R$-diagonal elements with non trivial kernel}\label{regard.bleu.en.haut.a.gauche.1} 
In the following, we shall use {\it polar decomposition} of non invertible elements 
of von Neumann algebras  (for example, in the following section, non invertible matrices). Recall that the polar decomposition of an element $x$ of a von Neumann algebra consists in writing $x=uh$, where $h\geq 0$ \st $\ker h=\ker x$, and $u$ is a partial isometry with initial space the orthogonal of $\ker x$ and with final space the closure of the image of $x$ (see the appendix of  \cite{dixmier} or the section 0.1 of \cite{sunder}). 

$R$-diagonal elements have been introduced by  Nica and Speicher in \cite{NS97}. In this section, we consider a $W^*$-noncommutative \pro space $(\mc{M},\tau)$.  In 1.9 of \cite{NS97}, $R$-diagonal elements of $(\mc{M},\tau)$ were characterized as the elements $x$ which can be written $x=uh$, where $u$ is a Haar unitary (i.e. $u$ is unitary, and for all $n\in \mathbb{Z}-\{0\}$, $\tau(u^n)=0$), and $h$ is a positive element $\tau$-free with $u$. If $x\in\mc{M}$ is $R$-diagonal and  if $x$ has a null kernel, then with the previous notations, $uh$ is the polar decomposition of $x$.   In the case where $x$ has a non trivial kernel, the polar decomposition of $x$ is $(up) h$, where $p$ is the projection on the orthogonal of $\ker (x)$. In this section, we shall prove that $up$, $h$ are free with amalgamation over the algebra $\Span\{p, 1-p\}$.

We first have to prove a preliminary result:  
\begin{propo}\label{regard.bleu.en.haut.a.gauche} Consider the space $(\A, \vfi)$ introduced in section \ref{amis+frere.26.07.05}, suppose moreover that  $(\A, \vfi)$ is a $W^*$-\pro space. Consider, in $\A$, a family $(y(s))_{s\in \n}$ of normal elements. 
Consider also a subalgebra $\B$ of $\A$ which contains $\D$ \st the family $(\B, (\{y(s)\})_{s\in\n})$ is $\vfi$-free. Then    the family $(\B, (\{y(s)\})_{s\in\n})$ is also free with amalgamation over $\D$.
\end{propo}

\begin{pr} 
Let $(\mc{N}, \tilde{\tau})$ be a $W^*$-\pro space  which is generated, as a $W^*$-algebra, by a family $x(s)$ ($s\in \n$) of $\tau$-semicircular elements and  an algebra $\widetilde{\B}$ isomorphic to $\B$ by a map $x\to \widetilde{x}$, \st the family $(\widetilde{\B}, (\{x(s)\})_{s\in\n})$ is $\tilde{\tau}$-free. The distributions of the $x(s)$'s are nonatomic, so  for each $s\in \n$, there exists a Borel function $f_s$ on the real line \st $f_s(x(s))$ has the same distribution as $y(s)$.   
Note that the $\tilde{\tau}$-freeness (resp. freeness with amalgamation over $\tilde{\D}$) of a family $(\A_i)_{i\in I}$ of $*$-subalgebras of $\mc{N}$ (resp. of  $*$-subalgebras of $\mc{N}$ which contain $\tilde{\D}$) is equivalent to the   $\vfi$-freeness (resp. the freeness with amalgamation over $\tilde{\D}$) of the family   $(\A_i'')_{i\in I}$ of von Neumann algebras they generate. So, by theorem  \ref{3.1.05.1} and by the fact that for all $s$, $f_s(x(s))\in \{x(s)\}''$, the family $(\tilde{\B}, (\{f_s(x(s))\})_{s\in\n})$ is free with amalgamation over $\D$. But the map \begin{eqnarray*}\B\cup \{y(s)\ste s\in \n\}& \to& \widetilde{\B}\cup \{f_s(x(s))\ste s\in \n\} \\ 
z&\mapsto &\begin{cases}\widetilde{z}& \textrm{if $z\in \B$}\\  f_s(x(s))& \textrm{if $z=y(s)$}
\end{cases} \end{eqnarray*} extends clearly to a $W^*$-\pro spaces isomorphism, hence the family $(\B, (\{y(s)\})_{s\in\n})$ is also free with amalgamation over $\D$.\end{pr}


\begin{cor}Consider  a $W^*$-noncommutative \pro space $(\mc{M},\tau)$, and an $R$-diagonal element $x$ of $\mc{M}$ with non trivial kernel. Let $p_1$ be the projection on $\ker x$, and $p_2=1-p_1$. Then the polar decomposition  $x=vh$ of $x$ is \st 
\begin{itemize}\item $v,h$ are free with amalgamation over $\D:=\Span\{p_1,p_2\}$,
\item $v$ has the $\D$-distribution of $up_2$, where $u$ is a Haar unitary $\tau$-free with $p_2$,
\item the $\D$-distribution of $h$ is defined by the fact that $h^2=x^*x$ and $p_2hp_2=h$.
\end{itemize}
Moreover, the projection on the final subspace of $v$ is $\tau$-free with $h$.
\end{cor}
 Note that this result could also have been deduced from lemma 2.6 of \cite{Sh00}, but the proof of this lemma is uncomplete, and a complete proof of the lemma would take as long as what we use to prove this corollary.

\begin{rmq}Elements with $\tau$-distributions such as the one of $v$ are called an {\it $(\alpha, \alpha)$-Haar partial isometries} in Remark 1.9 $3^\circ$ of \cite{NSS01}
\end{rmq}

\begin{pr}
By 1.9 of \cite{NS97}, $x$ can be written $x=uh$, where $u$ is a Haar unitary $\tau$-free with $h$, and the polar decomposition of $x$ is  $(up_2)h$, where $p_2$ is the projection on the orthogonal of $\ker x=\ker h$. Thus, with the notation $\D:=\Span\{1-p_2,p_2\}$, it suffices to prove the freeness with amalgamation over $\D$ of $up_2$ and $h$, which follows from the  freeness with amalgamation over $\D$ of $u$ and $h$, which follows from proposition \ref{regard.bleu.en.haut.a.gauche}. The projection on the final subspace of $up_2$ is $up_2u^*$ is $\tau$-free with $h$ by lemma 3.7 of \cite{haag2}.
\end{pr}

\begin{rmq} In the same way, we can prove the following: let $q_1$ be the projector on $\overline{\mathrm{Ran}(x)}$ and $q_2=1-q_1$, then the polar decomposition of $x$ is $wh$, where \begin{itemize}\item $w,h$ are free with amalgamation over $\D:=\Span\{q_1,q_2\}$,
\item $w$ has the $\D$-distribution of $q_1u$, where $u$ is a Haar unitary $\tau$-free with $q_1$,
\item the $\D$-distribution of $h$ is defined by the fact that $h^2=x^*x$ and $p_2hp_2=h$.
\end{itemize}
Moreover, the projection on $\ker w$ is $\tau$-free with $\D$.
\end{rmq}

\section{Asymptotic  freeness with amalgamation over $\D$ of rectangular random matrices}\label{oslo.mon.amour}
Since in the present section, we will prove asymptotic freeness with amalgamation over $\D$ of random matrices in an analogous way to the proofs of \cite{voic91} and \cite{voic98}, we shall frequently refer to those papers.  

Consider, for $n\geq 1$, $q_1(n)$,..., $q_d(n)$ positive integers with sum $n$ such that $\f{q_1(n)}{n}\ninf \rho_1$,..., $\f{q_d(n)}{n}\ninf \rho_d$ (recall that $\rho_1=\vfi(p_1)$,..., $\rho_d=\vfi(p_d)$). Then for all $n$, $\D$ can be identified with a $*$-subalgebra of the algebra $\Mn$ of complex $n\tii n$ matrices by $$\forall \la_1,\ldots,\la_d\in\C,\quad\diag(\la_1,\ldots,\la_d)\simeq\begin{bmatrix}\la_1I_{q_1(n)}& & \\  &\ddots & \\ & & \la_dI_{q_d(n)}.\end{bmatrix}$$The image of each $p_k$ will be denoted by $p_k(n)$. 
 $\tr$ will denote the normalized trace on $\Mn$, while $\Tr$ will denote the trace. $e(i,j;n)$ will denote the matrix-units of $\Mn$. 
\\
We shall refer to $M_n$ as a set of $n\tii n$ random matrices (over a \pro space not mentioned here), while the elements of $\Mn$ (which is a subalgebra of $M_n$) will be called constant matrices. $M_n$ is endowed with the state $\E(\tr(.))$, and the identification of $\D$ with a $*$-subalgebra of $M_n$ allows us to speak of convergence in $\D$-distribution of random matrices. 

The following result is an immediate corollary of Theorem 2.2 of \cite{voic98} and of 
theorem \ref{3.1.05.1}.
\begin{Th}\label{asfaogue}Let, for $s\geq 0$, $n\geq 1$, $Y(s,n)=\sum_{1\leq i,j\leq n}a(i,j;n,s)e(i,j;n)$ be a random matrix. Assume that $a(i,j;n,s)=\overline{a(j,i;n,s)}$ and that $$\{\Re a(i,j;n,s)\ste 1\leq i\leq j\leq n, s\in\n \}\cup\{\Im a(i,j;n,s)\ste 1\leq i< j\leq n, s\in\n\}$$ are independent Gaussian random variables, which are $(0,(2n)^{-1})$ if $i<j$ and $(0,n^{-1})$ if $i=j$. Let further $(B(j,n))_{j\in\n}$ be a family of elements of $\Mn$, stable under multiplication and adjonction, which contains $p_1(n)$,..., $p_d(n)$, such that for all $j$, the sequence $(||\ED(B(j,n)||)_n$ is bounded, and which converges in $\D$-distribution. 

Then the family $(\{B(j,n)\ste j\in\n\}, (\{Y(s,n)\})_{s})$ is asymptotically \fao as $n\to\infty$, and the limit $\D$-distribution of each $Y(s,n)$ is the $\D$-semicircular $\D$-distribution with covariance $\vfi$.  
\end{Th}

 \subsubsection*{Comparison with Theorem 4.1 of \cite{Sh96}} Since the $B(j,n)$'s are diagonal and since  convergence in $\D$-distribution is  less restricitve than convergence for $||.||_\infty$, which is the one used by  Shlyakhtenko in \cite{Sh96}, this  result cannot be deduced from Theorem 4.1 of \cite{Sh96}.

In order to modelize asymptotic $\D$-distribution of non hermitian gaussian random matrices, let us introduce the {\it $\D$-circular distribution with covariance $\vfi$}. It is the  $\D$-distribution of an element of $\A$ which can be written $c=\f{a+ia'}{\sqrt{2}}$, with $a,a'$ $\D$-semicircular
elements with covariance $\vfi$, which are \faop Note that the $\D$-distribution of $c$ can be defined by the following rules:
\begin{itemize}\item[(i)] $c_1(c)=c_1(c^*)=0$,
\item[(ii)] $\forall d\in\D,c_2(cd\otimes c)=c_2(c^*d\otimes c^*)=0, c_2(c d\otimes c^*)=c_2(c^*d\otimes c)=\vfi(d)$,  
\item[(iii)] $\forall k\geq 3, \eps_1,\ldots,\eps_k\in\{.,*\}, d_1,\ldots,d_{k}\in\D, c_k(c^{\eps_1}d_1\otimes\cdots \otimes c^{\eps_k}d_k)=0$.
\end{itemize}

\begin{cor}\label{1.1.05.88} The hypothesis are the same as the one of the previous theorem, except that the random matrices are not self-adjoint anymore, and their law is defined by the fact that $$\{\Re a(i,j;n,s)\ste 1\leq i,j\leq n, s\in\n \}\cup\{\Im a(i,j;n,s)\ste 1\leq i, j\leq n, s\in\n\}$$ are independant gaussian random variables, which are $(0,(2n)^{-1})$. Then the family $(\{B(j,n)\ste j\in\n\}, (\{Y(s,n)\})_{s})$ is asymptotically \fao as $n\to\infty$, and the limit $\D$-distribution of each $Y(s,n)$ is the $\D$-circular distribution with covariance $\vfi$.  
\end{cor}

\begin{pr} It suffices to notice that if $Y,Y'$ are independent random matrices as in the hypothesis of the previous theorem, then $\f{Y+iY'}{\sqrt{2}}$ has the distribution of the ones of  the hypothesis of the corollary.
\end{pr}

The previous corollary allows us to modelize asymptotic collective behaviour of independent rectangular gaussian random matrices with different sizes: consider, for   $s\geq 0$, $k,l\in[d]$, $n\geq 1$, $M(s,k,l,n)$ a random matrix of size $q_k(n)\tii q_l(n)$, with independent complex gaussian entries. In order to have a non trivial limit for asymptotic singular values of the $M(s,k,l,n)$'s (the {\it singular values} of a $q\tii q'$ matrix $M$ are the eigenvalues of of $MM^*$ if $q\leq q'$, and of $M^*M$ if $q\geq q'$), it is well known, by results about  Wishart matrices (see, e.g., \cite{hiai},\cite{pasturlejay}) that the variance of the entries must have the order of $q_k(n)$, i.e. of $n$. So we will suppose that the real and imaginary parts of the entries of the $M(s,k,l,n)$'s are  independent $N(0,(2n)^{-1})$. To give the asymptotic behaviour of these matrices amounts to give the asymptotic normalized traces of words of the type: 
\begin{equation}\label{tel.carmel}M(s_1,k_1,l_1,n)^{\eps_1}M(s_2,k_2,l_2,n)^{\eps_2}\cdots M(s_m,k_m,l_m,n)^{\eps_m},\end{equation}
where $m\geq 1$, $s_1,\ldots,s_m\in\n$, $\eps_1,\ldots,\eps_m\in\{.,*\}$, $k_1,l_1,\ldots,k_m,l_m\in[d]$ \st the product is possible and gives a square matrix.

In order to avoid problems of definition of the products, let us embed all this matrices in $n\tii n$ matrices: for all $(s,k,l,n)$, $M(s,k,l,n)$ will be replaced by $X(s,k,l,n):=p_k(n)Y(s,k,l,n)p_l(n)$, where $Y(s,k,l,n)$ is a random matrix as in the hypothesis of the previous corollary. Then if the product (\ref{tel.carmel}) is not defined, the product  \begin{equation}\label{tel.carmel.prime}X(s_1,k_1,l_1,n)^{\eps_1}X(s_2,k_2,l_2,n)^{\eps_2}\cdots X(s_m,k_m,l_m,n)^{\eps_m}\end{equation} is zero. In the other case,  the product (\ref{tel.carmel.prime}) is a simple element of $M_n$ (simple refers to the definition given in section \ref{amis+frere.26.07.05}), whose only non zero block is (\ref{tel.carmel}). If moreover, (\ref{tel.carmel}) is a square matrix, its normalized trace is the only non zero coordinate of $$\ED\lf[\lf(X(s_1,k_1,l_1,n)^{\eps_1}X(s_2,k_2,l_2,n)^{\eps_2}\cdots X(s_m,k_m,l_m,n)^{\eps_m}\ri)\ri].$$
  
So the  following corollary gives an answer to the question of the asymptotic collective behavior of independent rectangular Gaussian random matrices with different sizes.
\begin{cor} Let, for $s\geq 0$, $k,l\in[d]$, $n\geq 1$, $$X(s,k,l,n)=p_k(n)\lf[\sum_{1\leq i,j\leq n}a(i,j;n,k,l,s)e(i,j;n)\ri]p_l(n)$$ be a random matrix. Assume that $$\{\Re a(i,j;n,k,l,s)\ste i,j\in[n],k,l\in [d], s\in\n \}\cup\{\Im a(i,j;n,k,l,s)\ste i, j\in[n],k,l \in [d],s\in\n\}$$ are independent Gaussian random variables, which are $(0,(2n)^{-1})$. Let further $(B(j,n))_{j\in\n}$ be a family elements of $\Mn$, which satisfies the same assumptions as in the hypothesis of Theorem \ref{asfaogue}.  

Then the family $(\{B(j,n)\ste j\in\n\}, (\{X(s,k,l,n)\})_{s,k,l})$ is asymptotically \fao as $n\to\infty$.  
\end{cor}

\begin{pr}It is an immediate consequence of the previous corollary and of the fact that freeness with amalgamation over $\D$ is preserved by multiplication by elements of $\D$.  
\end{pr}


For $n\geq 1$, the set of matrices  $U$ of $p_k(n)\Mn p_k(n)$ \st $UU^*=U^*U=p_k(n)$ will be denoted by $\mathbb{U}_k(n)$. It is a compact group, isometric to the group of $q_k(n)\tii q_k(n)$ unitary matrices. By lemma  4.3.10 p. 160 of \cite{hiai}),  the partial isometry of the polar decomposition of $X(s,k,k,n)$ is {\it uniform} on $\mathbb{U}_k(n)$ (i.e. distributed according to the Haar measure).

\begin{propo}\label{1.01.05.1} Let, for $n\geq 1$, $V(s,k,n)$ ($s\in\n,k\in [d]$), be a family of independent random matrices, \st for all $s,k$, $V(s,k,n)$ is uniform on  $\mathbb{U}_k(n)$. Let  further $(B(j,n))_{j\in\n}$ be a family elements of $\Mn$ which satisfies the same assumptions as in the previous results. Then the family $ (\{B(j,n)\ste j\in\n\}, (\{V(s,k,n)\})_{s,k})$ is asymptotically \fao as $n\to\infty$.  
\end{propo}

Note that two elements of respectively $\A_{kk}$, $\A_{ll}$, with $k\neq l$ (or more generally of $\A_{kk'}$, $\A_{ll'}$, with $\{k,k'\}\cap\{l,l'\}=\emptyset$) are always free with amalgamation over $\D$, and that elements of $\A_{kk}$ are free with amalgamation over $\D$ \ssi they are free in the compressed space $(\A_{kk}, \vfi_k)$. So without the set of constant matrices, Proposition \ref{1.01.05.1} would be an easy consequence of Theorem 3.8 of \cite{voic91}. That being said, the proof of This Proposition is very closed to the one of  Theorem 3.8 of \cite{voic91}. 

The proof of the proposition relies on the following lemma. We endow $M_n$ with the norms $|M|_p=\lf(\E\tr(MM^*)^{r/2}\ri)^{\ff{r}}$. They  fulfill H\"older inequalities (see \cite{Nelson}).
\begin{lem} Let, for $n\geq 1$, $(M(i,n))_{i\in I}$ be a family of $n\tii n$ random matrices, and $(B(j,n))_{j\in\n}$ be a family elements of $\Mn$, which satisfies the same assumptions as in the hypothesis of Theorem \ref{asfaogue}.  Suppose moreover that for all $i\in I, r\geq 1$, the sequence $|M(i,n)|_r$ is bounded.
Suppose that for all $\delta>0$ and $n\geq 1$, their exists a family $(M(i,n,\delta))_{i\in I}$ of random $n\tii n$ matrices \st 
\begin{itemize}
\item[(i)] the family  $ (\{B(j,n)\ste j\in\n\}, (\{M(i,n,\delta)\})_{i\in I})$ is asymptotically \fao as $n\to\infty$,  
\item[(ii)] for all $i\in I, r\geq 1$, $\ds\overline{\lim}_{n\to \infty}|M(i,n,\delta)- M(i,n)|_r:= C(i,r,\delta)$ is such that $C(i,r,\delta)\underset{\delta\to 0}{\longrightarrow}0$.
\end{itemize}
Then the family  $ (\{B(j,n)\ste j\in\n\}, (\{M(i,n)\})_{i\in I})$ is asymptotically \fao as $n\to\infty$.  
\end{lem}
\begin{pr}Note first that if a sequence $(D_n)$ in $\D$ is \st for all positive $\delta$, there is a sequence $(D_n(\delta))$ in $\D$ which converges and \st $\ds\overline{\lim}_{n\to \infty}||D_n(\delta)- D_n||:= C(\delta)$ tends to zero as $\delta$ tends to zero, then $(D_n)$ is Cauchy, and hence converges.
So, by H\"older inequalities, the family  $ ((B(j,n)_{j\in\n}, (M(i,n))_{i\in I})$ has limit $\D$-distribution as $n\to\infty$.   Moreover, H\"older inequalities  implie too that the later limit $\D$-distribution is the limit, for convergence in $\D$-distribution, as $\delta$ tends to zero, of the limit $\D$-distribution of $ ((B(j,n))_{j\in\n}, (M(i,n,\delta))_{i\in I})$ as $n\to\infty$.  But the set of $\D$-distributions  of families $ ((b_j)_{j\in\n}, (m_i)_{i\in I})$  \st the family  $(\{b_j\ste j\in\n\}, (\{m_i\})_{i\in I})$ is \fao is obviously closed, so the lemma is proved. 
\end{pr}

Let us now give the proof of proposition \ref{1.01.05.1}.

\begin{pr} Consider independent random matrices $X(s,k,n)$ ($s\in\n,k\in [d],n\in\n$), \st for all $s,k,n$,  $X(s,k,n)$ has the same distribution as $X(s,k,k,n)$ of the previous corollary. Then, as noted before, one can suppose that for all $s,k,n$,  $V(s,k,n)$ is the partial isometry of the polar decomposition of $X(s,k,k,n)$. 
\\
In this proof, we shall use a particular fonctionnal calculus with the matrices  $X(s,k,n)^*X(s,k,n)$. Note that this matrices are simple elements of $M_n$ (simple refers to the definition given in section \ref{amis+frere.26.07.05}): they belong to $p_k(n)M_n p_k(n)$.  Here we shall ``erase'' the action of a function $f$ on the orthogonal of the image of the projector $p_k(n)$. This means that  $f(X(s,k,n)^*X(s,k,n))$ will mean  $p_k(n)f(X(s,k,n)^*X(s,k,n))p_k(n)$. So we can write: $$\ds V(s,k,n)=\lim_{\eps\to 0^+}X(s,k,n)(\eps+X(s,k,n)^*X(s,k,n))^{-1/\!2}.$$

{\it Step I. } As stated in the Step I of the proof of Theorem 3.8 of \cite{voic91}, there exists $C>0$ \st for all continuous bounded function $f: [0,\infty)\to (0,\infty)$, for all polynomial $P$, all $s,k$, and all $r\geq 1$, we have $$\ds\overline{\underset{n\to \infty}{\lim}}|X(s,k,n)P(X(s,k,n)^*X(s,k,n))-X(s,k,n)f(X(s,k,n)^*X(s,k,n)) |_r\leq C\!\!\!\sup_{0\leq t\leq C}\!\!\!|P(t)-f(t)|.$$
 
{\it Step II. } Consider $\eps >0$, and let $Y(s,k,n,\eps)=X(s,k,n)(\eps+X(s,k,n)^*X(s,k,n))^{-1/\! 2}$. We claim that  the family $$ (\{B(j,n)\ste j\in\n\}, (\{Y(s,k,n,\eps)\})_{s,k})$$ is asymptotically \fao as $n\to\infty$. It is an easy application of the lemma, using, for all positive $\delta$, the random matrices $$X(s,k,n)P_\delta(X(s,k,n)^*X(s,k,n)),$$ where $P_\delta$ is a polynomial \st $\sup_{0\leq t\leq C}|P_\delta(t)-(\eps+t)^{-1/\! 2}|\leq \delta$. Let us prove that the hypothesis of the lemma are satisfied. For $s,k,r$, the boundness of  the sequence $(|Y(s,k,n,\eps)|_r)_n$ comes from the boundness of the function $t\mapsto (t/(\eps+t))^{1/\! 2}$ on the positive half line, (i) is due to the previous corollary, and (ii) follows from step I.

{\it Step III. }  The conclusion is another application of the lemma, where $Y(s,k,n,\eps)$'s will play the roll of $M(i,n,\delta)$'s (and $\eps$ the roll of $\delta$). Let us, again, prove that the hypothesis of the lemma are satisfied. For $s,k,r$, the sequence $(|V(s,k,n)|_r)_n$ is bounded because matrices are in  $\mathbb{U}_k(n)$, and (i) follows from step II. Let us prove (ii). We have \begin{eqnarray*}|Y(s,k,n,\eps)-V(s,k,n)|_r&=&   
|V(s,k,n)(X(s,k,n)^*X(s,k,n))^{1/\!2}(\eps+X(s,k,n)^*X(s,k,n))^{-1/\!2}\\ &&\quad -V(s,k,n)|  
\\ 
&=&|(X(s,k,n)^*X(s,k,n))^{1/\!2}(\eps+X(s,k,n)^*X(s,k,n))^{-1/\!2}-p_k(n))|_r\end{eqnarray*}
$$\textrm{But }\quad \quad 0\leq t\leq \eps^{1/\!2}\qquad\Rightarrow\qquad  0\leq t^{1/\!2}/(\eps+t)^{1/\!2}\leq 1, $$
$$\textrm{and }\quad \quad t\geq \eps^{1/\!2}\qquad\Rightarrow \qquad 1-\f{\eps^{1/\!2}}{2(1+\eps^{1/\!2})}\leq t^{1/\!2}/(\eps+t)^{1/\!2} \leq 1.$$ So, if $Q(n)$ denotes the spectral projection of $X(s,k,n)^*X(s,k,n)$ for $[0, \eps^{1/\!2}]$,  \begin{eqnarray*}|V(s,k,n)-Y(s,k,n,\eps)|_r&\leq & \E(\tr (p_k(n)Q(n)p_k(n)))+\f{\eps^{1/\!2}}{2(1+\eps^{1/\!2})},\end{eqnarray*}
so (ii) is checked, since it is known (see \cite{szarek} for a precise result) that there exists a constant $C'$ \st $$\ds \overline{\lim}_{n\to \infty}\E(\tr (p_k(n)Q(n)p_k(n)))\leq C'\eps^{1/\!2}.$$
\end{pr}

Note that, as to Remark 2.3 of \cite{voic98}, it is obvious that the previous proposition also holds for subsequences of the natural numbers.  It allows us to prove the following corollary. Its proof is along the same lines as the one of Corollary 2.6 of \cite{voic98}: it relies on the fact that the topology of convergence in $\D$-distribution, for denombrable families, has denombrable bases of neighborhoods, and hence if for all $n$, $X(i,n)$ ($i\in I$ denombrable) is a family of random matrices with norms uniformly bounded, one can extract a subsequence $k(1)<k(2)<\ldots$ \st the family $(X(i,k(n)))_{i\in I}$ converges in $\D$ distribution.

\begin{cor} Let, for $n\geq 1$, $V(s,k,n)$ ($s\in\n,k\in [d]$), be a family of independent random matrices, \st for all $s,k$, $V(s,k,n)$ is uniform on  $\mathbb{U}_k(n)$ and let $F(\n\tii [d])\ni g\mapsto V^g(n)\in M_n$ be the semigroup morphism which sends the $(s,k)$-th generator into $V(s,k,n)$. Then, given $N\in\n,R>0$,and $g_0,\ldots,g_{N}\in F(\n\tii [d])-\{e\}$, the quantity 
$$\ds\sup\{|\tr(V^{g_0}(n)B_1V^{g_1}(n)B_2\cdots V^{g_{N-1}}(n)B_NV^{g_{N}}(n))|\ste \forall k\in [N],B_k\in \Mn ,||B_k||\leq R, \ED( B_k)=0\}$$ tends to zero as $n\to\infty$.
\end{cor}

At last, in the same way, one can translate the proof of Theorem 2.7 of \cite{voic98} to prove the following proposition:
\begin{propo}\label{bierbaum} Let, for $n\geq 1$, $V(s,k,n)$ ($s\in\n,k\in [d]$), be a family of independent random matrices, \st for all $s,k$, $V(s,k,n)$ is uniform on  $\mathbb{U}_k(n)$ and let $F(\n\tii [d])\ni g\mapsto V^g(n)\in M_n$ be the semigroup morphism which sends the $(s,k)$-th generator into $V(s,k,n)$. Fix $N\in\n$ and $R>0$. Let, for each $n\in\n$, $B_1(n),\ldots,B_N(n)$ be $n\tii n$ constant matrices \st for all $k\in [N]$, $||B_k(n)||\leq R$ and  $\ED(B_k(n))=0$. Then, given   $g_0,\ldots,g_N\in F(\n\tii [d])-\{e\}$ and $\eps>0$, the probability of the event $$\{|| \ED(V^{g_0}(n)B_1V^{g_1}(n)B_2\cdots V^{g_{N-1}}(n)B_NV^{g_{N}}(n))||\leq \eps\}$$ tends to $1$ as $n$ goes to infinity.
\end{propo}


\section{Analogue of free entropy of simple elements: the microstates approach}
\subsection{Definitions} For $q,q'$ positive integers, we will denote by $\mf{M}_{qq'}$ ($\mf{M}_q$ when $q=q'$) the set of $q\tii q'$ complex matrices.   

From now on, we suppose $(\A,\vfi)$ to be a tracial $W^*$-\pro space, endowed with a family of projectors such as presented in section \ref{amis+frere.26.07.05}: $\A$ is endowed with a family $p_1,\ldots, p_d$ of self-adjoint pairwise orthogonal projectors with sum $1$, and for all $k\in [d]$, $\rho_k$ denotes $\vfi(p_k)$. If $a$ is a simple element of $\A$, the unique $(k,l)\in [d]^2$ \st $a\in \A_{k,l}$ will be called the {\it type} of $a$. 

In this section, we shall define the entropy of families of simple elements of $\A$ as the asymptotic logarithm of the Lebesgue measure of sets of matrices with closed noncommutative moments. Thus we have to define, for all $n$, a set of matrices $\mf{M}_n(a)$ where we shall choose the microstates associated to a simple element $a$. We let, for $k\in [d]$, $q_k(n)$ be the integer part of $\rho_k n$. With this definition of $q_1(n)$, \ldots, $q_d(n)$, we keep the notations introduced in the beginning  of section \ref{oslo.mon.amour}. We define, for 
$(k,l)\in [d]^2$ and $n$ positive integer, $$\mf{M}_n(k,l)=p_k(n)\Mn p_l(n).$$ $\mf{M}_n(k,l)$ is endowed with the Lebesgue measure arising from the Euclidean structure defined by $<M,N>=\Re(\Tr M^*N).$ The norm arising from this Euclidean structure will be denoted by $||.||_2$, whereas $||.||$ still denotes the operator norm associated to the canonical hermitian norm on $\C^n$. 
We denote in the same time, without distinction, by $\La$ the tensor product of this measures on any product of such spaces.

Consider $a_1,\ldots a_N\in \A$ simple elements with respective types $(k(1), l(1)),\ldots, (k(N),l(N))$. Let us define, for  $n, r$ positive integers, $\eps, R$ positive numbers, $$\Gamma_R(a_1,\ldots, a_N; n, r, \eps)$$ the set of families $(A_1,\ldots, A_N)\in \mf{M}_n(k(1), l(1))\tii\cdots \tii \mf{M}_n(k(N),l(N))$ \st  for all $i=1,\ldots,N$, $||A_i||\leq R$ and for all $p\in \{1,\ldots, r\}$, for all $i_1,\ldots, i_p\in [N]$, for all $\eps_1,\ldots \eps_p\in \{*,.\}$, $$||\ED(A_{i_1}^{\eps_1}\cdots A_{i_p}^{\eps_p})-\ED(a_{i_1}^{\eps_1}\cdots a_{i_p}^{\eps_p})||\leq \eps.$$ 
 Let us  then define $$\chi^\D_R(a_1,\ldots, a_N; r, \eps)=\ds\limsup_{n\to\infty}\ff{n^2}\log \La\lf(\Gamma_R(a_1,\ldots, a_N; n, r, \eps)\ri)+L\log n+D,$$ where \begin{eqnarray*}L&=&\ds
\sum_{k,l=1}^d\rho_k\rho_l\lf|\{i\in [N] \ste a_i\textrm{ has type $(k,l)$}\}\ri|,\\ D&=& \ds\sum_{k=1}^d\rho_k^2\log\rho_k 
\lf|\{i\in [N] \ste a_i\textrm{ has type $(k,k)$}\}\ri|   .\end{eqnarray*}
We define then  $$\chi^\D_R(a_1,\ldots, a_N)=\ds\inf_{r,\eps}\chi^\D_R(a_1,\ldots, a_N; r, \eps),$$ and at last, $$\chi^\D(a_1,\ldots, a_N)=\ds\sup_R\chi^\D_R(a_1,\ldots, a_N).$$

\subsection{Particular cases, comparison with already defined quantities} 

\subsubsection{Case where all $a_i$'s are of the same type $(k,k)$}Then $$\chi^\D_R(a_1,\ldots, a_N; r, \eps)=\rho_k^2\chi_R^\mathrm{Voic}(a_1,\ldots, a_N; r, \eps/\rho_k),$$where $\chi_R^\mathrm{Voic}(a_1,\ldots, a_N; r, \eps)$ stands for $\chi_R(a_1,\ldots, a_N; r, \eps)$ as it is defined   in section 1.2 of \cite{voic98} or p. 279 of \cite{hiai}, when considering $a_1,\ldots, a_N$ as non-selfadjoint elements of $(\A_{kk}, \vfi_k)$.

\subsubsection{Case where $N=1$ and $a_1$ has type $(k,l)$, with $k\neq l$, $q_k(n)=q_l(n)$ for all $n$}\label{malade.mais.battant.puis.bien}
Then $\Gamma_R(a; n, r, \eps)$ is the set of matrices $A\in \Mn (k,l)\simeq \mf{M}_{q_k(n)}$ \st $||M||\leq R$ and for all  $s$ positive integer \st $2s\leq r$, the $s$-th moment of the spectral law of $AA^*$ is within $\eps$ with the $s$-th moment of the distribution of $aa^*$ in $(\A_{kk},\vfi_k)$.
Thus $\chi^\D(a)=\rho_k^2(\chi_+^\mathrm{Voic}(aa^*)-\log \rho_k)$, where $\chi^\mathrm{Voic}_+(aa^*)$ is defined  p. 282 of \cite{hiai}, when considering $aa^*$ as a positive element of $(\A_{kk}, \vfi_k)$.

Note that, if $\mu$ is the distribution of $aa^*$ in $(\A_{kk}, \vfi_k)$, then $\chi^\D_+(aa^*)=\int\!\!\!\int\log |x-y|\ud \mu(x)\ud\mu(y)+\log \pi +3/\!2$.

\subsubsection{Case where $N=1$ and $a_1$ has type $(k,l)$, with $k\neq l$}\label{16.2.05.5}
Then $\Gamma_R(a; n, r, \eps)$ is the set of matrices $A\in \Mn (k,l)$ \st $||M||\leq R$ and for all  $s$ positive integer \st $2s\leq r$, the $s$-th moment of the spectral law of $AA^*$ is within $\eps$ with the $s$-th moment of the distribution of $aa^*$ in $(\A_{kk},\vfi_k)$ and    the $s$-th moment of the spectral law of $A^*A$ is within $\eps$ with the $s$-th moment of the distribution of $a^*a$ in $(\A_{ll},\vfi_l)$.  Let us  define $$\mu=\begin{cases} \textrm{distribution of $aa^*$ in $(\A_{kk},\vfi_k)$}&\textrm{if $\rho_k\leq \rho_l$},\\  \textrm{distribution of $a^*a$ in $(\A_{ll},\vfi_l)$}&\textrm{if $\rho_k>\rho_l$}. \end{cases}$$ Then we prove the following proposition (proof 
postponed in the appendix): 

\begin{propo}\label{livre.cardio.pauline} If $R^2$ is more than the supremum  of the support of $\mu$, \begin{equation}\label{14.02.05.1}\chi^\D_R(a)=\alpha^2\Sigma(\mu)+(\beta-\alpha)\alpha\!\!\int\log x\ud \mu(x)+\alpha\beta(\log \f{\pi}{\alpha}+1)+\f{\alpha^2}{4}-\alpha^2\!\!\int_{\f{\beta-\alpha}{\alpha}}^{\f{\beta}{\alpha}}x\log x\ud x,\end{equation} where $\Sigma(\mu)=\int\!\!\!\int\log |x-y|\ud \mu(x)\ud\mu(y)$, $\alpha=\min\{\rho_k,\rho_l\}$ and $\beta=\max\{\rho_k,\rho_l\}$.\end{propo}
 
We can verify that when $\alpha=\beta$, this formula coincides with the formula of $\chi^\D(a)$ given by \ref{malade.mais.battant.puis.bien}.


\subsection{Preliminary lemma}The following lemma is a very useful tool for integration on sets of rectangular matrices. It gives the ``law'' of the singular values of rectangular matrices distributed according to the Lebesgue measure. Its proof 
is postponed in the appendix.  Consider $1\leq q\leq q'$ integers. Denote by $\mf{M}_{q,q'}$ the set of $q\tii q'$
 complex matrices. Denote by $\mc{U}_q$  the group of $q\tii q$  unitary matrices, and by $T_q$ the torus of diagonal matrices of $\mc{U}_q$. As an homogeneous space, $\mc{U}_q/T_q$ is endowed with a unique distribution invariant under the left action of $\mc{U}_q$, denoted by  $\gamma_q$. 
Denote by $\mc{U}_{q,q'}$ the set of matrices $v$ of $\mf{M}_{q,q'}$ which satisfy $vv^*=I_q$ (i.e. whose lines are orthogonal with norm $1$).  As an homogeneous space (under the right action of $\mc{U}_{q'}$), $\mc{U}_{q,q'}$ is endowed with a unique distribution invariant under this actions, denoted by $\gamma_{q,q'}$. Note that   $\gamma_{q,q'}$ is also invariant under the left action of $\mc{U}_q$.
\begin{lem}\label{16.2.05.2} 
Define $$\R^q_{+,<}=\{x\in \R^q\ste 0<x_1<\cdots <x_q\}.$$  Then the map \begin{eqnarray*}\Psi\, :\,\mc{U}_q/T_q\,\tii\, \R^q_{+,<}\,\tii\, \mc{U}_{q,q'}&\to& \mf{M}_{q,q'} \\
(uT_q,x,v)&\mapsto & u\diag (x_1,\ldots, x_q)^{1/\! 2}u^*v\end{eqnarray*}is injective onto a set with complement of null Lebesgue measure. Moreover,  the push-forward, by $\Psi^{-1}$, of the Lebesgue measure, is $\gamma_q\otimes \sigma_{q,q'}\otimes  \gamma_{q,q'}$, where $\sigma_{q,q'}$ is the \pro measure on  $\R^q_{+,<}$ with density \begin{equation}\label{banque.avec.david}\ds\f{\pi^{qq'}}{\prod_{j=1}^{q-1}j!\prod_{j=q'-q}^{q'-1}j!}\Delta(x)^2  \prod_{i=1}^{q}x_j^{q'-q}.\end{equation}
\end{lem}

{\bf Remark.} It will be more useful to use the following consequence of the lemma. Let us denote, for $G$ compact group, $\Haar(G)$ the Haar  \pro measure on $G$. The measures $\gamma_q$ and $\gamma_{q,q'}$ are push-forwards of $\Haar(\mc{U}_q)$, $\Haar(\mc{U}_{q'})$ by the respective maps $u\to uT_q$, $v\to Pv$, where $P$ is the $q\tii q'$ matrix with diagonal entries equal to $1$, and other entries equal to $0$. Then the map  \begin{eqnarray*}\tilde{\Psi}\, :\,\mc{U}_q\,\tii\, (\R^+)^q\,\tii\, \mc{U}_{q'}&\to& \mf{M}_{q,q'} \\
(u,x,v)&\mapsto & u\diag (x_1,\ldots, x_q)^{1/\! 2}u^*Pv\end{eqnarray*}
is surjective and the push-forward of the measure $\Haar (\mc{U}_q)\otimes \ff{q!}\tilde{\sigma}_{q,q'}\otimes \Haar (\mc{U}_{q'})$
by $\tilde{\Psi}$ is the Lebesgue measure, where $\tilde{\sigma}_{q,q'}$ is the measure on $(\R^+)^n$ with density given by formula (\ref{banque.avec.david}). 

\subsection{Classical properties of entropy}
The following properties are analogous to properties of Voiculescu's entropy, the proofs are analogous too, and the straightforward adaptation will be left to the reader. For the proof of proposition \ref{16.2.05.1}, the classical change of variables formula used for square matrices needs to be replaced by the result given below  lemma \ref{16.2.05.2}.

\begin{propo}$\chi^\D$ is subadditive: for  $1\leq m <N$, $$\chi^\D(a_1,\ldots a_N)\leq \chi^\D(a_1,\ldots a_M)+\chi^\D(a_{M+1},\ldots a_N).$$
\end{propo}

\begin{propo}\label{16.2.05.1}For $R_1>R>\max\{||a_1||,\ldots, ||a_N||\}$, we have $$\chi^\D_{R_1}(a_1,\ldots a_N)=\chi^\D_R(a_1,\ldots a_N).$$
\end{propo}

\begin{propo}[Upper semicontinuity]\label{15.2.05.1}
Consider, for $m\geq 1$, $a_{m,1},\ldots, a_{m,N}$ simple elements of $\A$ \st \begin{itemize}\item for all $i$, $a_{m,i}$ has the same type as $a_i$,
 \item the family $(a_{m,1},\ldots, a_{m,N})$ converges  in $\D$-distribution to $(a_1,\ldots, a_N)$,
\item for all $i$, the sequence $||a_{m,i}||$ is bounded. 
\end{itemize}
Then $\chi^\D(a_1,\ldots, a_N)\geq \limsup \chi^\D(a_{m,1},\ldots, a_{m,N}).$
\end{propo} 
For example, the proposition holds 
if  for all $i$, $a_{m,i}$ has the same type as $a_i$ and converges strongly to $a_i$. 

\begin{propo}
Consider $b_1,\ldots, b_N$ simple elements \st for all $i\in [N]$, $b_i$ has the type of $a_i$ and $b_i-a_i\in \{a_1,\ldots, a_{i-1}\}''$. Then $\chi^\D(a_1,\ldots, a_N)=\chi^\D(b_1,\ldots, b_N)$.  
\end{propo}

\subsection{Entropy and freeness with amalgamation over $\D$}
Adaptating the section 5 of \cite{entropyII} and using proposition \ref{bierbaum}, we obtain the following result:
\begin{Th}\label{17.2.05.1} If the simple elements $a_1, \ldots, a_N$ are free with amalgamation over $\D$, then $$\chi^\D(a_1,\ldots, a_N)=\chi^\D(a_1)+\cdots +\chi^\D(a_N).$$ 
\end{Th}

\subsection{Change of variable formula} Consider $a_1,\ldots a_N\in \A$ simple elements with respective types $(k(1), l(1)),\ldots, (k(N),l(N))$. In this section, for $i\in [N]$, in order to simplify expressions, we denote $\mf{M}_n(k(i),l(i))$ by $\mf{M}_n(a_i)$.  Since we are going to work with adjoints of the $a_i$'s, we have to define, for $i\in [N]$, $\eps\in\{*,1\}$, $(k(i,\eps), l(i,\eps))$ to be the type of $a_i^\eps$ i.e. $$   (k(i,\eps), l(i,\eps))=\begin{cases}(l(i),k(i))&\textrm{if $\eps=*$}\\ (k(i),l(i))&\textrm{if $\eps=1$}\end{cases}$$

Define $\mc{F}$ to be the set of formal power series in the noncommutative variables $X_1,X_1^*,\ldots, X_N,X_N^*$ endowed with the natural involution $F\to F^*$. Let us define, for $m\geq 0, i_1,\ldots, i_m\in [N],\eps_1,\ldots, \eps_m\in\{1,*\}$, $C_{\substack{ i_1,\ldots, i_m\\ \eps_1,\ldots,\eps_m}}$ the map from $\mc{F}$ to $\C$ which maps a series $F$ to its coefficient in $X_{i_1}^{\eps_1}\cdots X_{i_m}^{\eps_m} $. 

A {\it multi-radius of convergence} for $F\in \mc{F}$ 
is a family $(R_1,\ldots, R_N)$ of positive numbers \st $$M(F; R_1,\ldots, R_N):=\sum_{m\geq 0}\sum_{\substack{i_1,\ldots, i_m\in [N]\\ \eps_1,\ldots,\eps_m\in\{*,1\}}} |C_{\substack{ i_1,\ldots, i_m\\ \eps_1,\ldots,\eps_m}}(F)|R_{i_1}\cdots R_{i_m}<\infty  .$$

Define, for 
 $i\in [N]$, the set $\mc{F}_{i}$ of formal power series $F\in \mc{F}$ \st for all $m\geq 1$, for all $i_1,\ldots, i_m , \eps_1,\ldots,\eps_m$, 
$$C_{\substack{ i_1,\ldots, i_m\\ \eps_1,\ldots,\eps_m}}(F)\neq 0\Rightarrow k(i)=k(i_1,\eps_1),l(i_1,\eps_1)=k(i_2,\eps_2),\ldots, l(i_{m-1},\eps_{m-1})=k(i_m,\eps_m),l(i_m,\eps_m)=l(i).$$

Consider $F=(F^{(1)},\ldots, F^{(N)})$ in $\mc{F}_1\tii \cdots \tii \mc{F}_N$. We suppose moreover that there exists $(R_1,\ldots, R_N)$ common multiradius of convergence of the $F^{(i)}$'s \st for all $i$, $||a_i||<R_i$.

Let also, for $n\geq 1$,  $F$ be  the map defined on the Cartesian product, for $i\in [N]$, of the open ball of $(M_n(a_i), ||.||)$ with center zero   and radius $R_i$, by \begin{eqnarray*}F\,:\,\ds\prod_{i=1}^NB_{M_n(a_i)}(0,R_i) &\to & \prod_{i=1}^NB_{M_n(a_i)}(0,R'_i)\\
(A_1,\ldots, A_N)&\mapsto& F( A_1,\ldots, A_N)\end{eqnarray*}

Then $F$ is analytic, and with the natural identification between the set $\mc{L}\lf(\ds\prod_{i=1}^NM_n(a_i)\ri)$ of endomorphisms of $\ds\prod_{i=1}^NM_n(a_i)$ and the Cartesian product $\ds\prod_{1\leq i,j\leq N} \mc{L}\lf(M_n(a_j), M_n(a_i)\ri)$, the differential $DF(A)$ of $F$ at $A=( A_1,\ldots, A_N)\in \ds\prod_{i=1}^NB_{M_n(a_i)}(0,R_i)$ has $(i,j)$-th block
$$\ds\sum_{\substack{m\geq 1\\ i_1,\ldots, i_m\in [N]\\ \eps_1,\ldots,\eps_m\in\{*,1\}}} C_{\substack{ i_1,\ldots, i_m\\ \eps_1,\ldots,\eps_m}}(F^{(i)})\bigg[\sum_{\substack{l\in [m]\\ i_l=j\\ \eps_l=1}} L( A_{i_1}^{\eps_1}\cdots A_{i_{l-1}}^{\eps_{l-1}})\circ R(A_{i_{l+1}}^{\eps_{l+1}}\cdots A_{i_{m}}^{\eps_{m}})$$ $$+\sum_{\substack{l\in [m]\\ i_l=j\\ \eps_l=*}}  L(A_{i_1}^{\eps_1}\cdots A_{i_{l-1}}^{\eps_{l-1}})\circ R(A_{i_{l+1}}^{\eps_{l+1}}\cdots A_{i_{m}}^{\eps_{m}})\circ \Adj\bigg],$$ where for $A$ matrix, $L(A)$ (resp. $R(A)$) denotes the operator of left (resp. right) multiplication by $A$, and $\Adj$ denotes the operator of adjonction in $\Mn$. 

Now, we are going to compute the Jacobian of $F$ at  $A=( A_1,\ldots, A_N)$. It is the absolute value of the determinant of the differential of $F$ at  $( A_1,\ldots, A_N)$, that is $$(\det DF(A)DF(A)^*)^{1/\! 2}=\exp\ff{2}\Tr \log DF(A)DF(A)^* ,$$ where the adjoint is taken when considering $DF(A)$ as an endomorphism of the space $\ds\prod_{i=1}^NM_n(a_i)$ endowed with the product euclidean structure. Note that the identification between $$\mc{L}\lf(\ds\prod_{i=1}^NM_n(a_i)\ri)\;\textrm{ and }\;\ds\prod_{1\leq i,j\leq N} \mc{L}\lf(M_n(a_j), M_n(a_i)\ri)$$ preserves the adjonction in the following way: 
$$\lf(\lf[M_{i,j}\ri]_{i,j=1}^N\ri)^*=\lf[M_{j,i}^*\ri]_{i,j=1}^N,$$ 
and composition in the following way 
$$\lf[L_{i,j}\ri]_{i,j=1}^N\circ \lf[M_{i,j}\ri]_{i,j=1}^N=\lf[\sum_{k=1}^n L_{i,k} \circ M_{k,j}\ri]_{i,j=1}^N.$$

Let, for $n\geq 1$,  $\mf{L}_n$ be the space of endomorphisms $L$ of  $\ds\prod_{i=1}^NM_n(a_i)$ such that, for all $i,j\in [N]$, $L_{i,j}$ is a linear combinaison of linear maps of the type $L_A\circ L_B$ (with $A\in \Mn (k(i),k(j))$ and $B\in   \Mn (l(j),l(i))$) and of maps of the type $L_A\circ L_B\circ \Adj$  (with $A\in \Mn (k(i),l(j))$ and $ B\in   \Mn (k(j),l(i))$).  $\mf{L}_n$ is a subalgebra of $\mc{L}(\ds\prod_{i=1}^NM_n(a_i))$ closed under adjonction. Indeed, we have, for all $X,Y,Z,T$ matrices with suitable sizes, \begin{eqnarray*}
(L(X)\circ R(Y))\circ (L(Z)\circ R(T))&=&L(XZ)\circ R(YT),\\ (L(X)\circ R(Y))\circ (L(Z)\circ R(T)\circ \Adj)&=&L(XZ)\circ R(YT)\circ \Adj,\\
(L(X)\circ R(Y)\circ \Adj) \circ (L(Z)\circ R(T))&=&L(XT^*)\circ R(YZ^*)\circ \Adj, \\ (L(X)\circ R(Y)\circ \Adj) \circ (L(Z)\circ R(T)\circ \Adj)&=&L(XT^*)\circ R(YZ^*),\\
 (L(X)\circ R(Y))^*&=&L(X^*)\circ R(Y^*),\\ (L(X)\circ R(Y)\circ \Adj )^*&=&L(Y)\circ R(X)\circ \Adj.\end{eqnarray*}
Thus in order to compute the Jacobian of $F$, it suffices to be able to compute the trace of a self-adjoint element of $\mf{L}_n$. Note that  the identification between $$\mc{L}\lf(\ds\prod_{i=1}^NM_n(a_i)\ri)\;\textrm{ and }\;  \ds\prod_{1\leq i,j\leq N} \mc{L}\lf(M_n(a_j), M_n(a_i)\ri)$$ preserves the trace in the following way: 
$$\Tr\lf(\lf[M_{i,j}\ri]_{i,j=1}^N\ri)=\ds\sum_{i=1}^N\Tr M_{i,i}.$$
Moreover, if  $\lf[M_{i,j}\ri]_{i,j=1}^N\in\mf{L}_n $ is self-adjoint, then for all $i\in \{1,\ldots, N\}$, $M_{i,i}=M_{i,i}^*$, and one can write
$$M_{i,i}=\ds\sum_{\alpha} c_\alpha (L_{X_\alpha}\circ R_{Y_\alpha}+L_{X_\alpha^*}\circ R_{Y_\alpha^*})+\sum_{\beta} c_\beta (L_{Z_\beta}\circ R_{T_\beta}\circ \Adj+L_{T_\beta^*}\circ R_{Z_\beta^*}\circ \Adj),$$ where $\alpha, \beta$ run in disjoint finite sets, and for all $\alpha,\beta$, $c_\alpha,c_\beta$ are real,   $$X_\alpha\in p_{k(i)}(n)\Mn p_{k(i)}(n), \quad Y_\alpha\in p_{l(i)}(n)\Mn p_{l(i)}(n),\quad  Z_\beta,T_\beta\in M_n(a_i).$$
Thus the trace of  $M_{i,i}$ is $$\ds\sum_{\alpha} c_\alpha 2\Re (\underbrace{\Tr X_\alpha\Tr Y_\alpha+\Tr X_\alpha^*\Tr Y_\alpha^*}_{\in \R}),$$ $$\textrm{i.e. }\quad  \ds\sum_{\alpha} 2c_\alpha (\Tr X_\alpha\Tr Y_\alpha+\Tr X_\alpha^*\Tr Y_\alpha^*).$$

Thus, in order to compute the Jacobian of $F$ at $A=(A_1,\ldots, A_N)$, we have to introduce the following objects.  Let $\mf{S}_2=\{e,\tau\}$ be the group of permutations of the set $\{1,2\}$, and $\C[\mf{S}_2]=\C e\oplus \C\tau$ be its convolution algebra.  For $j\in [N]$, define the $\C$-linear map $$D_j\, :\, \mc{F}\, \to\, \mc{F}\otimes\mc{F}\otimes  \C[\mf{S}_2]$$   by   $$D_j(X_{i_1}^{\eps_1}\cdots X_{i_m}^{\eps_m})=\sum_{\substack{l\in [m]\\ i_l=j\\ \eps_l=1}}  X_{i_1}^{\eps_1}\cdots X_{i_{l-1}}^{\eps_{l-1}}\otimes X_{i_{l+1}}^{\eps_{l+1}}\cdots X_{i_{m}}^{\eps_{m}}\otimes e+\sum_{\substack{l\in [m]\\ i_l=j\\ \eps_l=*}}  X_{i_1}^{\eps_1}\cdots X_{i_{l-1}}^{\eps_{l-1}}\otimes X_{i_{l+1}}^{\eps_{l+1}}\cdots X_{i_{m}}^{\eps_{m}}\otimes \tau.$$
Define then, for any $*$-algebra $\mc{M}$ (which will be $\mc{F}$, or $\Mn$, or $\A$), 
the product and the adjonction on the linear space $ \mc{M}\otimes\mc{M}\otimes  \C[\mf{S}_2]$ defined by the rules: 
\begin{eqnarray*}
(X\otimes Y\otimes e)\tii (Z\otimes T\otimes e)&=&XZ\otimes YT\otimes e,\\ (X\otimes Y\otimes e)\tii ( Z\otimes T\otimes \tau)&=&XZ\otimes YT\otimes \tau,\\
(X\otimes Y\otimes \tau )\tii (Z\otimes T\otimes e)&=&XT^*\otimes YZ^*\otimes \tau, \\ (X\otimes Y\otimes \tau )\tii(Z\otimes T\otimes \tau)&=&XT^*\otimes YZ^*\otimes e,\\
 (X\otimes Y\otimes e)^*&=&X^*\otimes Y^*\otimes e,\\ (X\otimes Y\otimes \tau )^*&=&Y\otimes X\otimes \tau.\end{eqnarray*}
 If moreover, the algebra $\mc{M}$ is endowed with a linear functional $f$ (which will be $\Tr$ if $\mc{M}=\Mn$, and $\vfi$ if $\mc{M}=\A$), then we shall endow $ \mc{M}\otimes\mc{M}\otimes  \C[\mf{S}_2]$ with the linear functional $2 f\otimes f\otimes \delta_e$, where $\delta_e$ is the state on $ \C[\mf{S}_2]$ defined by $\delta_e(e)=1$, $\delta_e(\tau)=0$. 

With this notations, if one uses the identifications, for $X,Y\in \Mn$:  
$$L(X)\circ R(Y)\simeq X\otimes Y\otimes \delta_e \in \Mn\otimes \Mn\otimes  \C[\mf{S}_2], \quad L(X)\circ R(Y)\circ \Adj \simeq X\otimes Y\otimes \tau \in \Mn\otimes \Mn\otimes  \C[\mf{S}_2],$$and thus identifies $\mc{L}_n$ with a subset of the finite-dimensional algebra $\mf{M}_N\otimes (\Mn\otimes \Mn\otimes \C[\mf{S}_2])$,
the Jacobian of $F$ at $A=(A_1,\ldots, A_N)$ is $$\exp\ff{2}\ds\Tr\otimes  \Tr\otimes \Tr \otimes \delta_e \lf(\log (DF (A)DF^*(A))\ri).$$

With this tools, adaptating the proof of proposition 3.5 of \cite{entropyII}, we have the following proposition. $\Delta_{\Tr\otimes \vfi\otimes \vfi\otimes (2\delta_e)}$ denotes the Kadison-Fuglede determinant of the linear functional $\Tr\otimes \vfi\otimes \vfi\otimes (2\delta_e)$ on $\mf{M}_N\otimes \A\otimes \A\otimes \C[\mf{S}_2]$ (see section 2 of \cite{haag2} for a brief introduction to Kadison-Fuglede determinant).

\begin{propo}Consider $F=(F^{(1)},\ldots, F^{(N)}), G=(G^{(1)},\ldots, G^{(N)})$ which both belong to $\mc{F}_1\tii \cdots \tii \mc{F}_N$, and \st for all $i\in [n]$, $$G^{(i)}(F^{(1)},\ldots, F^{(N)})=X_i.$$We suppose moreover that there exists $(R_1,\ldots, R_N)$ common multiradius of convergence of the $F^{(i)}$'s and $(R'_1,\ldots, R'_N)$ common multiradius of convergence of the $G^{(i)}$'s \st \begin{itemize} \item[(i)] for all $i$, $||a_i||<R_i$, \item[(ii)] for all $i$, $M(F^{(i)}; R_1,\ldots,R_N)<R_i'$. \end{itemize}
Then $$\chi^\D(F^{(1)}(a_1,\ldots, a_N),\ldots, F^{(N)}(a_1,\ldots, a_N))$$ $$\geq \log \Delta_{\Tr\otimes \vfi\otimes \vfi\otimes (2\delta_e)}[D_jF^{(i)}(a_1,\ldots, a_N)]_{\substack{i,j=1}}^N+\chi^\D(a_1,\ldots,a_N).$$
If moreover, $$(M(G^{(1)}; R_1',\ldots, R_N'), \ldots , M(G^{(N)}; R_1',\ldots, R_N'))$$ is a common multiradius of the $F^{(i)}$'s, then we have equality.\end{propo}

\begin{rmq}Note that if instead of (ii) we have  $M(F^{(i)}; ||a_1||,\ldots,||a_N||)<R_i'$ for all $i$, then one can reduce the $R_i$'s in order to have (i) and (ii).
\end{rmq}

As a corollary, we have the following result, whose proof is an adaptation of the proof of proposition 6.3.3 of \cite{hiai}.

\begin{cor}\label{18.2.05.1} Let $P_1,\ldots, P_N$ be noncommutative polynomials of $X_1,X_1^*,\ldots, X_N^*$ \st for all $i\in [N]$, $a_i+P_i(a_1,\ldots,a_N)$ has type $(k,l)$. Then  for $\alpha$ sufficiently near $0$, $$\chi^\D(a_1+\alpha P_1(a_1,\ldots, a_N),\ldots, a_N+\alpha P_N(a_1,\ldots, a_N))$$
$$=\log \Delta_{\Tr\otimes \vfi\otimes \vfi\otimes (2\delta_e)}[D_jF_\alpha^{(i)}(a_1,\ldots, a_N)]_{\substack{i,j=1}}^N+\chi^\D(a_1,\ldots, a_N),$$where for all $i$, $F_\alpha^{(i)}=X_i+P(X_1,\ldots, X_N)$. \end{cor} 

\subsection{Functional calculus and entropy}\label{tu.me.touche.mais.je.ne.pense.pas.que.ce.soit.tres.reel}
For $x\in \A$, with spectral decomposition $uh$, and $f$ real Borel function on $[0,\infty)$, bounded on the spectrum of $h$, let $f(x)$ denote $uf(h)$. If $f(0)=0$ and $f$ is positive on $(0,\infty)$, then  the polar decomposition of $f(x)$ is $uf(h)$. 

We begin with the following lemma,      analogous to lemma 6.3.5 of \cite{hiai}.
\begin{lem}Consider $k,l\in [d]$ \st $k\leq l$,  $a\in \A_{kl}$ \st the distribution $\mu$ of $(aa^*)^{1/\! 2}$ in $(\A_{kk},\vfi_k)$ satisfies $$\Sigma(\mu)>-\infty, \quad \int\log t\ud\mu(t)>-\infty,$$ and $f$ a continuous increasing function on $[0,\infty)$, \st $f(0)=0$ and $f$ is positive on $(0,\infty)$. Then there exists a sequence $(f_m)$ of smooth  functions on $[0,\infty)$, \st for all $m$, $f_m(0)=0$, $f'_m$ is positive on $[0,\infty)$, $$||f_m(a)-f(a)||\ninf 0,\quad \chi^\D(f_m(a))\ninf \chi^\D(f(a)).$$
\end{lem}  

\begin{pr}
Let us prove that there exists a sequence $(f_m)$ of smooth  functions on $[0,\infty)$, \st for all $m$, $f_m(0)=0$, $f'_m$ is positive on $[0,\infty)$, $$||f_m(a)-f(a)||\ninf 0,\quad \chi^\D(f_m(a))\ninf \chi^\D(f(a)).$$
Then we will have, by upper semicontinuity: $\chi^\D(f(a))\geq \limsup \chi^\D(f_m(a))$, and by (\ref{14.02.05.1}), it will suffice to prove $$\Sigma(f(\mu))\leq \liminf\Sigma(f_m(\mu)),\quad \int\log(f(t))\ud \mu(t)\leq \liminf\int\log(f_m(t))\ud \mu(t),$$ where for all function $g$, $g(\mu)$ denotes the push-forward of $\mu$ by $g$.

Consider, for $m\geq 1$, $\delta(m)\in (0,1/\! m)$ \st $$\int\!\!\!\int_{|t-s|<\delta(m)}\log|t-s|\ud \mu(t)\ud \mu(s)\geq -\ff{m}, \quad \int_{[0,\delta(m)]}\log t\ud \mu(t)\geq -\ff{m}, $$
 $$\int\!\!\!\int_{|t-s|<\delta(m)}\ud \mu(t)\ud \mu(s)\leq \ff{m\log m}, \quad \mu([0,\delta(m)])\leq \ff{m\log m}. $$Let us extend $f$ by the value $0$ on the negative real numbers. Let, for $m\geq 1$, $\phi_m$ be a nonnegative smooth function with support in $[0, 1/\!m]$ \st $\int \phi_m(t)\ud t=1$ and $$|f*\phi_m(t)-f(t)|\leq \f{\delta(m)}{m}$$ for all $t\in [0,S]$, where $S$ is the maximum of the support of $\mu$. Define $f_m(t):=\f{t}{m}+f*\phi_m(t)$. Since $f$ is increasing and $\phi_m\geq 0$, $f*\phi_m$ is increasing, hence $f_m'\geq 1/\! m$. Moreover, $f_m(0)=0$ and $f_m$ converges uniformly to $f$ on $[0,S]$. 

Similarly to p. 267 of \cite{hiai}, we can prove that for $m$ large enough to satisfy \st 
$$\forall s,t \in [0, S], |t-s|<\delta(m)\Rightarrow |f(t)-f(s)|<1,$$   
we have $ \Sigma(f_m(\mu))\geq \Sigma(f(\mu))-\f{2}{m}.$ Moreover, 
for $t\in [0,\delta(m)]$, 
$f_m(t)\geq t/\! m$ and for 
$t\in [\delta(m),S]$, 
$$f_m(t)\geq \f{\delta(m)}{m}+f(t)-|f*\phi_m(t)-f(t)|\geq f(t),$$ so if $\delta(m)\leq 1$,  \begin{eqnarray*}\int \log (f_m(t))\ud \mu(t)&\geq & \int_{[\delta(m),S]}\log(f(t))\ud \mu(t)+\int_{[0,\delta(m)] } \log(T/\! m)\ud\mu(t)\\   &\geq & \int_{[0,S]}\log(f(t))\ud \mu(t)+\int_{[0,\delta(m)] } \log t\ud\mu(t)-\log(m)\mu([0,\delta(m)])\\ &\geq & \int \log(f(t))\ud \mu(t)-\f{2}{m},  \end{eqnarray*}what closes the proof. \end{pr}

Adaptating the proof of proposition 6.3.6 of \cite{hiai} (with the density of eigenvalues presented in lemma \ref{16.2.05.2}), we obtain the following proposition: 

\begin{propo}\label{18.02.05.6}Consider $f_1,\ldots, f_N$ continuous increasing functions on $[0,\infty)$, with value $0$ in $0$, and positive on $(0,\infty)$.  Consider $a_1,\ldots, a_N$ simple elements of $\A$  \st for all $i$, $\chi^\D(a_i)>-\infty$. Then $$\chi^\D(f_1(a_1)),\ldots, f_N(a_N))\geq \chi^\D(a_1,\ldots, a_N)+\sum_{i=1}^N \chi^\D(f_i(a_i))-\chi^\D(a_i).$$ Moreover, equality holds if the functions are strictly increasing.
\end{propo}

\subsection{Maximization of free entropy}
\subsubsection{One variable}\label{17.2.05.2}  The problem here is to maximize $\chi^\D(a)$, where $a$ is taken among a set of simple elements of type $(k,l)$ in $\A$. 
For such $a$, we have seen in \ref{16.2.05.5} a formula which express $\chi^\D(a)$ as a function of $\mu$ (distribution of $aa^*$ if $\rho_k\leq \rho_l$, distribution of $a^*a$ in the other case), of $\rho_k$  and of $\rho_l$. For example, let us suppose that $\rho_k\leq \rho_l$. 
\begin{propo}Fix $c$ positive. Then among elements $a$ of type $(k,l)$ \st $\vfi_k(aa^*)\leq c$, the maximizers of $\chi^\D$ are those for which the distribution of $aa^*$ in $(\A_{kk}, \vfi_k)$ is the push-forward, by $t\to \f{c\rho_k}{\rho_l}t$, of the Marchenko-Pastur distribution with parameter $\f{\rho_l}{\rho_k}$. 
\end{propo}

\begin{pr}First of all, since for all $a$, for all $\la>0$, $\chi^\D(\la a)=\chi^\D(a)+2\rho_k \rho_l\log \la$, it suffices to prove it when $c=\f{\rho_l}{\rho_k}$.
According to proposition \ref{livre.cardio.pauline},  $\chi^\D(a)$ is maximal if and only if $$A(\mu):=\Sigma(\mu)+\lf(\f{\rho_l}{\rho_k}-1\ri)\int\log x\ud \mu(x)$$ is maximal. Note that the condition $\vfi_k(aa^*)\leq c$ is equivalent to $\int x\ud \mu(x)\leq c$. But proposition 5.3.7 of \cite{hiai} states that the functional $$B(\mu):=\Sigma(\mu)+\lf(\f{\rho_l}{\rho_k}-1\ri)\int\log x\ud \mu(x)-\int x\ud \mu(x)$$ is maximized, among \pro measures on $[0,\infty)$, by the Marchenko-Pastur distribution $\mu_c$ with parameter $c$. So  for all $\mu$ \pro measure on $[0,\infty)$ \st  $\int x\ud \mu(x)\leq c$, \begin{eqnarray*}A(\mu_c)-A(\mu)&=&B(\mu_c)-B(\mu)+\int x\ud \mu_c(x)-\int x\ud \mu(x)\\  &=&B(\mu_c)-B(\mu)+(c-\int x\ud \mu(x))\\ &\geq & 0,   \end{eqnarray*} with equality \ssi $\mu=\mu_c$. 
\end{pr}

\subsubsection{$N$ variables} Fix $k\neq l\in [N]$ \st $\rho_k\leq \rho_l$.
A consequence of subadditivity, of theorem \ref{17.2.05.1}, and of  the previous section is the fact that if $c_1,\ldots, c_N$ are positive numbers, among $N$-tuples $(a_1,\ldots, a_N)$ of type $(k,l)$ elements which satisfy $$\forall i\in [N],\quad \vfi_k(a_ia_i^*)=c_i,$$ the maximum of $\chi^\D(a_1,\ldots, a_N)$ is realized on free with amalgamation over $\D$ families  $(a_1,\ldots, a_N)$ of elements \st  the distribution of each $a_ia_i^*$ in $(\A_{kk},\vfi_k)$ is  the push-forward, by $t\to \f{c_i\rho_k}{\rho_l}t$, of the Marchenko-Pastur distribution with parameter $\f{\rho_l}{\rho_k}$.  The following theorem states the reciprocal to this fact. 
\begin{Th}If the maximum is realized on a family $(a_1,\ldots, a_N)$, then the family is free with amalgamation over $\D$, and  the distribution of each $a_ia_i^*$ in $(\A_{kk},\vfi_k)$ is  the push-forward, by $t\to \f{c_i\rho_k}{\rho_l}t$, of the Marchenko-Pastur distribution with parameter $\f{\rho_l}{\rho_k}$.  
\end{Th}

\begin{pr} {\it Step I. } First of all, since for all $a_1,\ldots, a_N$ of type $(k,l)$, for all $\la_1,\ldots, \la_N>0$, \begin{equation}\label{regard.qui.perce.et.qui.met.quelque.chose.en.feu.dans.moi.c.est.un.peu.bebete}\chi^\D(\la_1 a_1,\ldots, \la_Na_N)=\chi^\D(a_1,\ldots, a_N)+2\rho_k \rho_l\log (\la_1\cdots\la_N),\end{equation} it suffices to prove it when each $c_i$ is $\f{\rho_l}{\rho_k}$.

{\it Step II. } By \teo \ref{regard.qui.perce.et.qui.met.quelque.chose.en.feu.dans.moi},  if $b$ is an element of type $(k,l)$ \st the distribution of $bb^*$ in $(\A_{kk},\vfi_k)$ is  the Marchenko-Pastur distribution with parameter $\f{\rho_l}{\rho_k}$, then for all $n$ even integer, $$c_n^{(k)}(b\otimes b^* \otimes \cdots \otimes b\otimes b^*)=\f{\rho_l}{\rho_k}\delta_{n,2}.$$
Let us compute the $\vfi$-distribution of a free with amalgamation over $\D$ family $(b_1,\ldots, b_N)$ of such elements. For all $i_1,\ldots, i_{2r+1}, i_0\in [N]$, we have \begin{eqnarray*}\vfi(b_{i_1}b_{i_2}^*\cdots b_{i_{2r+1}}b_{i_0}i^*)&=&   \rho_k\vfi_k(b_{i_1}b_{i_2}^*\cdots b_{i_{2k+1}}b_{i_0}^*)\\
&=&   \rho_k\sum_{\substack{0\leq j\leq r\\ i_{2j+1}=i_0}}\vfi_k(b_{i_1}\cdots b_{i_{2j}}^*)c_2^{(k)}(b_{i_0}\otimes \underbrace{\ED (b_{2j+2}^*\cdots b_{2r+1})}_{\vfi_l(b_{2j+2}^*\cdots b_{2r+1}).p_l} b_{i_0}^*)\\ 
&=&   \rho_k\sum_{\substack{0\leq j\leq r\\ i_{2j+1}=i_0}}\vfi_k(b_{i_1}\cdots b_{i_{2j}}^*)c_2^{(k)}(b_{i_0}\otimes \vfi_l(b_{2j+2}^*\cdots b_{2r+1}).\underbrace{p_l b_{i_0}^*}_{b_{i_0}^*})\\  
&=&   \rho_k\sum_{\substack{0\leq j\leq r\\ i_{2j+1}=i_0}}\vfi_k(b_{i_1}\cdots b_{i_{2j}}^*)c_2^{(k)}(b_{i_0}\otimes \vfi_l(b_{2j+2}^*\cdots b_{2r+1}).b_{i_0}^*)\\
&=& \rho_k\sum_{\substack{0\leq j\leq r\\ i_{2j+1}=i_0}}\vfi_k(b_{i_1}\cdots b_{i_{2j}}^*)\vfi_l(b_{2j+2}^*\cdots b_{2r+1})c_2^{(k)}(b_{i_0}\otimes b_{i_0}^*)\\
&=& \rho_l\sum_{\substack{0\leq j\leq r\\ i_{2j+1}=i_0}}\vfi_k(b_{i_1}\cdots b_{i_{2j}}^*)\vfi_l(b_{2j+2}^*\cdots b_{2r+1})
\end{eqnarray*}

{\it Step III. } Now, consider a family $(a_1,\ldots, a_N)$ of elements of type $(k,l)$ \st for all $i$, $\vfi_k(a_ia_i^*)=\f{\rho_l}{\rho_k}$ and the maximum of the entropy is realized on $(a_1,\ldots, a_N)$. Let us prove that the $\D$-distribution of the family is the one  of $(b_1,\ldots, b_N)$ of step II. Since the elements are simple and for all $x$, $\vfi(x^*)=\overline{\vfi(x)}$, it suffices to prove that for all $r\geq 0$, for all  $i_1,\ldots, i_{2r+1}, i_0\in [N]$, we have \begin{equation*}\vfi(a_{i_1}a_{i_2}^*\cdots a_{i_{2r+1}}a_{i_0}^*)=\rho_l\sum_{\substack{0\leq j\leq r\\ i_{2j+1}=i_0}}\vfi_k(a_{i_1}\cdots a_{i_{2j}}^*)\vfi_l(a_{2j+2}^*\cdots a_{2r+1}).\end{equation*}

Fix $\la\in \C$ and  let us define the polynomial $P=\la X_{i_1}X_{i_2}^*\cdots X_{i_{2r+1}}$, and $d=P(a_1,\ldots, a_N)$. We have, for $\alpha$ sufficiently near $0$, $$\chi^\D(a_1,\ldots, a_N)\geq \chi^\D(a_1,\ldots, a_{i-1}, \rho_l^{1/\! 2} \f{a_{i_0}+\alpha d}{(\rho_k\vfi_k((a_{i_0}+\alpha d)(a_{i_0}+\alpha d)^*))^{1/\! 2}},a_{i+1},\ldots, a_N).$$  Thus the derivative with respect to $\alpha$, at $\alpha=0$, of the difference between right-hand side and left-hand side  of the previous equation is zero. According to (\ref{regard.qui.perce.et.qui.met.quelque.chose.en.feu.dans.moi.c.est.un.peu.bebete}) and to corollary \ref{18.2.05.1}, for $\alpha$ sufficiently near $0$,  the difference between  right-hand side and left-hand side  of the previous equation is equal to $$\log \Delta_{\Tr\otimes \vfi\otimes \vfi\otimes (2\delta_e)}[D_jF_\alpha^{(i)}(a_1,\ldots, a_N)]_{\substack{i,j=1}}^N-\rho_k\rho_l\log\f{\rho_k\vfi_k((a_{i_0}+\alpha d)(a_{i_0}+\alpha d)^*)}{\rho_l} ,$$where for all $j$, $F_\alpha^{(j)}=X_j+\alpha\delta_{i,j}P$. We have$$  \log \Delta_{\Tr\otimes \vfi\otimes \vfi\otimes (2\delta_e)}[D_jF_\alpha^{(i)}(a_1,\ldots, a_N)]_{\substack{i,j=1}}^N$$ $$=\Tr\otimes \vfi\otimes \vfi\otimes \delta_e\log\lf([D_jF_\alpha^{(i)}(a_1,\ldots, a_N)]_{\substack{i,j=1}}^N[D_iF_\alpha^{(j)}(a_1,\ldots, a_N)^*]_{\substack{i,j=1}}^{N}\ri)$$
 $$=\Tr\otimes \vfi\otimes \vfi\otimes \delta_e\log\lf(I_N\otimes1\otimes 1\otimes e+\alpha(A+A^*)+\alpha^2AA^*\ri),$$ where $A$ is the $N\tii N$ matrix with $(i,j)$-th entry $$\begin{cases}0&\textrm{if $i\neq i_0$,}\\ D_jP(a_1,\ldots, a_N)&\textrm{if $i= i_0$.}
\end{cases} $$Thus  $$  \log \Delta_{\Tr\otimes \vfi\otimes \vfi\otimes (2\delta_e)}[D_jF_\alpha^{(i)}(a_1,\ldots, a_N)]_{\substack{i,j=1}}^N$$ $$= \alpha\vfi\otimes \vfi\otimes \delta_e (D_{i_0}P(a_1,\ldots, a_N)+(D_{i_0}P(a_1,\ldots, a_N))^*)+o(\alpha).$$ Thus $$\f{\partial}{\partial \alpha}_{|\alpha=0}\log \Delta_{\Tr\otimes \vfi\otimes \vfi\otimes (2\delta_e)}[D_jF_\alpha^{(i)}(a_1,\ldots, a_N)]_{\substack{i,j=1}}^N$$ $$=\vfi\otimes \vfi\otimes \delta_e (D_{i_0}P(a_1,\ldots, a_N)+(D_{i_0}P(a_1,\ldots, a_N))^*)=2 \Re\vfi\otimes \vfi\otimes \delta_e (D_{i_0}P(a_1,\ldots, a_N)) .$$ Moreover, since $\vfi_k(a_{i_0}a_{i_0}^*)= \f{\rho_l}{\rho_k}$, 
 $$\f{\partial}{\partial \alpha}_{|\alpha=0}\rho_k\rho_l\log\f{\rho_k\vfi_k((a_{i_0}+\alpha d)(a_{i_0}+\alpha d)^*)}{\rho_l} =\rho_k^2\vfi_k(a_{i_0}d^*+da_{i_0}^*)=2\rho_k^2\Re \vfi_k(P(a_1, \ldots, a_N)a_{i_0}^*).$$ Thus $\vfi\otimes \vfi\otimes \delta_e (D_{i_0}P(a_1,\ldots, a_N))$ and $ \rho_k^2\vfi_k(P(a_1, \ldots, a_N)a_{i_0}^*)$ have the same real part. Recall that   $P=\la X_{i_1}X_{i_2}^*\cdots X_{i_{2r+1}}$. What we did is true for any $\la\in \C$, so   $$ \rho_k^2\vfi_k(P(a_1, \ldots, a_N)a_{i_0}^*)=\vfi\otimes \vfi\otimes \delta_e (D_{i_0}P(a_1,\ldots, a_N)).$$
Now recall the definition of $D_{i_0}P$ and choose $\la=1$. This gives $$ \rho_k^2\vfi_k( a_{i_1}a_{i_2}^*\cdots a_{i_{2r+1}}a_{i_0}^*)
=\sum_{\substack{0\leq j\leq r\\ i_{2j+1}=i_0}}\vfi(a_{i_1}\cdots a_{i_{2j}}^*)\vfi(a_{2j+2}^*\cdots a_{2r+1}),$$i.e. $$ \rho_k\vfi( a_{i_1}a_{i_2}^*\cdots a_{i_{2r+1}}a_{i_0}^*)
=\sum_{\substack{0\leq j\leq r\\ i_{2j+1}=i_0}}\rho_k\vfi_k(a_{i_1}\cdots a_{i_{2j}}^*)\rho_l\vfi_l(a_{2j+2}^*\cdots a_{2r+1}),$$which closes the proof.
\end{pr}

\begin{cor}Consider a family $(a_1,\ldots, a_N)$ of elements of type $(k,l)$ \st $$\chi^\D(a_1,\ldots, a_N)=\chi^\D(a_1)+\cdots+\chi^\D(a_N)>-\infty.$$ Then the family is free with amalgamation over $\D$.\end{cor}

\begin{pr}We use the notation defined in the beginning of section \ref{tu.me.touche.mais.je.ne.pense.pas.que.ce.soit.tres.reel}. Since for all $i$, $\chi^\D(a_i)>-\infty$, the distribution of $a_ia_i^*$ in $(\A_{kk},\vfi_k)$ is nonatomic, hence the distribution of $|a_i|$ in  in $(\A_{kk},\vfi_k)$ is also nonatomic. Hence we can find a continuous increasing function $f_i$ on $[0,\infty)$, with value $0$ in $0$, and positive on $(0,\infty)$, \st the distribution of  $f_i(a_i)f_i(a_i)^*$ is the Marchenko-Pastur distribution with parameter $\f{\rho_l}{\rho_k}$.  Then by proposition \ref{18.02.05.6} and subadditivity, $\chi^\D(f_1(a_1),\ldots, f_N(a_N))=\chi^\D(f_1(a_1))+\cdots +\chi^\D(f_N(a_N)) $, hence by the previous theorem, $f_1(a_1),\ldots, f_N(a_N)$ are free with amalgamation over $\D$, and so do $a_1,\ldots, a_N$, because for all $i$, $a_i\in \{f_i(a_i)\}''$. 
\end{pr}

\subsubsection*{Question} It would be interesting to have a characterization of $R$-diagonal elements with non trivial kernel involving the entropy defined in this paper. Inspired by the papers \cite{HP99}, \cite{NSS99.2}, we ask the following question: given a compactly supported \pro measure $\nu$ on $\R^+$, what are the elements $a\in \A p_1$ \st $aa^*$ has distribution $\nu$ in $(\A_{22},p_2)$, and \st  $\chi^\D(p_1ap_{2},p_2ap_{2}, p_2a^*p_{1}, p_2a^*p_{2})$ is maximal ?

\section{Analogue of free Fisher's information for simple elements: the microstate-free approach}

In this
section, we present a notion of free Fisher's information for simple elements, constructed without using the microstates, like what was done by Voiculescu in \cite{entropyV} and by Shlyakhtenko in \cite{Sh00} for elements of a $W^*$-\pro space. For a synthetic presentation of the free Fisher's information of  elements of a $W^*$-\pro space, see section 2 of \cite{NSS99}.

\subsection{Definitions}
In this section, we suppose that $(\A,\vfi)$ is a $W^*$-\pro space, with $\vfi$ faithful tracial state. $L^2(\A)$ will denote the Hilbert space obtained by completing $\A$ for the norm $||a||_2=(\vfi(aa^*))^{\ff{2}}$, $a\in \A$. 

$\A$ acts on $L^2(\A)$ on the right and on the left, so one can define, for $k,l \in [d]$, $ L^2(\A)_{kl}=p_k L^2(\A) p_l.$ We still have an identification between   $ L^2(\A)$ and $$\begin{bmatrix}L^2(\A)_{11}&\cdots & L^2(\A)_{1d}\\ \vdots& &\vdots\\ L^2(\A)_{d1}&\cdots & L^2(\A)_{dd}\end{bmatrix},$$
by the map $x\in L^2(\A)\to \begin{bmatrix}p_1 x p_1 &\cdots & p_1 xp_d \\ \vdots & & \vdots \\ p_d x p_{1} & \cdots & p_dx p_{d}\end{bmatrix}.$ We still call the non zero elements of $\ds \underset{k,l\in[d]}{\cup} L^2(\A)_{kl}$ the {\it simple elements} of $L^2(\A)$. We define, for $a$ non zero simple element of $L^2(\A)$, $r(a)$ ($r$ is for {\it row}) and $c(a)$ ($c$ is for {\it column}) the unique numbers of $[d]$ \st $$a\in L^2(\A)_{r(a), c(a)}.$$

Moreover, $x\to x^*$ extends to $L^2(\A)$, and for all $a,b\in \A, x\in L^2(\A)$, we have $(axb)^*=b^*x^*a^*$. For all $k\in [d]$, the state $\vfi_k$ on $\A_{kk}$ extends to $
 L^2(\A)_{kk}$, so the conditional expectation $\ED$ extends to $L^2(\A)$, and we still have $$\forall d,d' \in \D, \forall x\in L^2(\A), \ED(dxd')=d\ED(x)d'.$$ In the same way, for $n\geq 1$ and $\pi\in \NC(n)$, we can extend $\ED_n$, $\ED_\pi$, $c_n$ and $c_\pi$ to $L^2(\A)\otimes (\A^{\otimes n-1}),$ and the relations (\ref{def.cum.1.1.05.1}), (\ref{def.cum.1.1.05.2}), (\ref{def.cum.1.1.05.3}), and (\ref{def.cum.1.1.05.4}) remain true. 

A family $(a_i)_{i\in I}$ of elements of $L^2(\A)$ is said to be a {\it self-adjoint family} if there exists an involution $*$ of $I$ \st for all $i\in I$, $a_i^*=a_{*(i)}$. 

A finite sequence $(a_1,\ldots, a_n)$ of simple elements of $L^2(\A)$ is said to be a {\it square sequence} if for all $i\in [n]$, $c(a_i)=l(a_{i+1})$ (with $a_{n+1}:=a_1$). 

\begin{Def}\label{Chao+Wozniak=Choose.the.Bosniaque} Let $(a_i)_{i\in I}$ be a self-adjoint family of simple elements of $\A$. A family $(\xi_i)_{i\in I}$ of simple elements of $L^{2}(\A)$ is said to {\it fulfill conjugate relations} for $(a_i)_{i\in I}$ if for all $i\in I$, $$r(\xi_i)=c(a_i),\quad\quad c(\xi_i)=r(a_i),$$ and if one of the following equivalent proposition is true : 
\begin{itemize}
\item[(i)] for all $n\geq 0$, for all $i,i_1,\ldots, i_n\in I$ \st $(\xi_i, a_{i_1},\ldots, a_{i_n})$ is a square sequence, $$\vfi_{r(\xi_i)}(\xi_i a_{i_1}\cdots a_{i_n})=\sum_{m=1}^n\delta_{i,i_m}\vfi_{c(\xi_i)}(a_{i_1}\cdots a_{i_{m-1}})\vfi_{r(\xi_i)}(a_{i_{m+1}}\cdots a_{i_n}),$$
\item[(ii)] \begin{itemize}\item[-] for all $i\in I$, $\ED(\xi_i)=0$,
\item[-] for all $i,j\in I$, $\ED(\xi_ia_j)=\delta_{ij}p_{r(\xi_i)}$, 
\item[-] for all $n\geq 2$, for all $i,i_1,\ldots, i_n\in I$ \st $(\xi_i, a_{i_1},\ldots, a_{i_n})$ is a square sequence, $$\ED (\xi_i a_{i_1}\cdots a_{i_n})=\sum_{m=1}^n\delta_{i,i_m}\eta_{r(\xi_i),c(\xi_i)}\circ\ED(a_{i_1}\cdots a_{i_{m-1}})\ED(a_{i_{m+1}}\cdots a_{i_n}),$$ (we recall that for all $k,l\in[d]$, $\eta_{k,l}$ is the involution of $\D$ which permutes the $k$-th and the $l$-th columns in the representation of elements of $\D$ as $d\tii d$ diagonal matrices matrices),
\end{itemize} 
\item[(iii)] for all $n\geq 0$, for all $i,i_1,\ldots, i_n\in I$ \st $(\xi_i, a_{i_1},\ldots, a_{i_n})$ is a square sequence, $$c_{n+1}(\xi_i\otimes a_{i_1}\otimes \cdots \otimes a_{i_n})=\delta_{n,1}\delta_{i,i_1}p_{r(\xi_i)}.$$ 
\end{itemize}  \end{Def}

\begin{pr}Consider a family $(\xi_i)_{i\in I}$ of simple elements of $L^{2}(\A)$ \st  for all $i\in I$, $$l_i:=r(\xi_i)=c(a_i),\quad\quad k_i:= c(\xi_i)=r(a_i),$$ and let us prove the equivalence of (i), (ii), and (iii). The equivalence between (i) and (ii) is obvious by definition of $\ED$: for all $x\in L^2(\A)$, $$\ED( x)=\ds\sum_{k=1}^d\vfi_k(p_kap_k).p_k.$$

Now, suppose (ii) true, and let us prove (iii) by induction on $n$. Note that $c_1=\ED_1$, and that for all $\xi\in L^2(\A)$, $a\in \A$, $$c_2(\xi\otimes a)=\ED(\xi a)-\ED(\xi)\ED(a),$$ 
so (iii) is proved for $n=0,1$. 
Now, suppose it proved for to all ranks $0,1,\ldots, n-1$, with $n\geq 2$, and let us prove it to the rank $n$. Consider $i,i_1,\ldots, i_n\in I$ \st $(\xi_i, a_{i_1},\ldots, a_{i_n})$ is a square sequence. We have 
\begin{eqnarray*}c_{n+1}(\xi_i\otimes a_{i_1}\otimes \cdots \otimes a_{i_n})&=&
\ED(\xi_i a_{i_1}\cdots  a_{i_n})-\ds\sum_\pi c_\pi(\xi_i\otimes a_{i_1}\otimes \cdots \otimes a_{i_n}),\end{eqnarray*}
where the sum is taken over all noncrossing partitions $\pi$ of $\{0,\ldots, n\}$ which are $<1_\{0,\ldots, n\}$, and in which all blocks are associated to square subsequences of $(\xi_i, a_{i_1},\ldots, a_{i_n})$. 

Consider such a partition $\pi$, and  apply the factorization formula (\ref{volanniv.2.3.05.1}) to $c_\pi(\xi_i\otimes a_{i_1}\otimes \cdots \otimes a_{i_n})$. If it is not null, then the block of $\pi$ containing $0$ has only one other element, say $m$, and in this case, we have $$c_\pi(\xi_i\otimes a_{i_1}\otimes \cdots \otimes a_{i_n})=\ds\!\!\!\!\prod_{\substack{V\in \pi_1\\ V=\{t_1<\cdots <t_r\} }}\!\!\!\!\eta_{l_i,l_{t_r}}\circ c_{r}(a_{i_{t_1}}\otimes \cdots \otimes a_{i_{t_r}})\!\!\!\!\!\!\prod_{\substack{V\in \pi_2\\ V=\{t_1<\cdots <t_r\} }}\!\!\!\!\!\!\eta_{l_i,l_{t_r}}\circ c_{r}(a_{i_{t_1}}\otimes \cdots \otimes a_{i_{t_r}}),$$ where $\pi_1$ (resp. $\pi_2$) is the partition induced by $\pi$ on  $\{1,\ldots, m-1\}$ (resp. $\{m+1,\ldots , n\}$), i.e. to $$\delta_{i,i_m}\eta_{l_i,k_i}\circ c_{\pi_1}(a_{i_1}\otimes\cdots \otimes a_{i_{m-1}})c_{\pi_2}(a_{i_{m+1}}\otimes\cdots\otimes a_{i_n}),$$
Thus we  have \begin{eqnarray*}c_{n+1}(\xi_i\otimes a_{i_1}\otimes \cdots \otimes a_{i_n})
&=&
\ED(\xi_i a_{i_1}\cdots  a_{i_n})-
\ds\sum_{m=1}^n\!\!\!\!\!\!\sum_{\substack{\pi_1\in \NC(m-1)\\ \pi_2\in\NC(\{m+1,\ldots, n\})}}  \\
&&\quad\quad\delta_{i,i_m}
\eta_{l_i,k_i}\circ c_{\pi_1}(a_{i_1}\otimes\cdots \otimes a_{i_{m-1}})c_{\pi_2}(a_{i_{m+1}}\otimes\cdots\otimes a_{i_n})\\ 
&=&
\ED(\xi_i a_{i_1}\cdots  a_{i_n})-\ds\sum_{m=1}^n\delta_{i,i_m}\eta_{l_i,k_i}\circ\ED(a_{i_1}\cdots a_{i_{m-1}})\ED(a_{i_{m+1}}\cdots a_{i_n})\\  &=&0.
\end{eqnarray*}

The reciprocal implication is analoguous.
\end{pr}

\begin{Def}Let $(a_i)_{i\in I}$ be a self-adjoint family of simple elements of $\A$.  
\\
\\
1. A family $(\xi_i)_{i\in I}$ of simple elements of $L^{2}(\A)$ is said to de a {\it conjugate system} for $(a_i)_{i\in I}$ if  it fulfills conjugate relations an if in addition we have that \begin{equation}\label{28.02.05.1} \forall i\in I,\quad\quad \xi_i\in \overline{\Alg (\{a_i\ste i\in I\}\cup \D)}^{||.||_2}\subset L^2(\A).\end{equation} 
\\
\\
2. Let the $\D$-Fisher's information of $(a_i)_{i\in I}$ be $$\Phi_r((a_i)_{i\in I})=\begin{cases}\ds\sum_{i\in I}||\xi_i||^2&\textrm{if $(\xi_i)_{i\in I}$ is a conjugate system,}\\ \infty&\textrm{if there is no conjugate system.}\end{cases} $$
\end{Def}

\subsubsection*{Remark}\label{retard.1.3.05.1} 
1. The algebra generated by $\{a_i\ste i\in I\}\cup \D$ is the set of elements of $\A$ which have a $d\tii d$ matrix representation of the type 
$$\begin{bmatrix}P_{11}(a_i\ste i\in I)& \cdots &P_{1d}(a_i\ste i\in I)\\\vdots&&\cdots\\  P_{d1}(a_i\ste i\in I)&\cdots & P_{dd}(a_i\ste i\in I) \end{bmatrix},$$where for all $k,l\in [d]$, $P_{kl}$ is a polynomial in the noncommutative variables $X_i\ste i\in I$, and $ P_{kl}(a_i\ste i\in I)\in \A_{kl}$. So the conjugate relations can be viewed as a prescription for the inner products in $L^2(\A)$ between $\xi_i$ ($i\in I$) and elements of this subalgebra. It follows that the conjugate system for $(a_i)_{i\in I}$ is unique, if it exists. Note moreover that the existence of the conjugate system is equivalent to the existence of any family in $L^2(\A)$ which fulfills the conjugate relations; indeed, if $(\xi_i)_{i\in I}$ fulfill the conjugate relations and if we set, for all $i\in I$,  $\gamma_i$ to be the projection of $\xi_i$ onto $\overline{\Alg (\{a_i\ste i\in I\}\cup \D)}^{||.||_2}$, then   $(\gamma_i)_{i\in I}$ will also fulfill the conjugate relations, hence will give a conjugate system.  

2. Consider an involution $*$ of $I$ \st for all $i\in I$, $a_{*(i)}=a_i^*$. Consider a family  $(\xi_i)_{i\in I}$ which fulfills conjugate relations. Then define, for all $i\in I$, $$\tilde{\xi}_i=\f{\rho_{c(a_i)}}{\rho_{r(a_i)}}\xi_{*(i)}^*.$$ Then   $(\tilde{\xi}_i)_{i\in I}$ fulfills conjugate relations, hence we have $\xi_i=\tilde{\xi}_i$ for all $i$ if  $(\xi_i)_{i\in I}$ is a conjugate system. It can be written \begin{equation}\label{dudududadadaisallIwanttosaytoyou}\xi_{*(i)}=\f{\rho_{r(a_i)}}{\rho_{c(a_i)}}\xi_i^*.\end{equation}

\begin{pr} Let us prove 2. It suffices to prove that $(\tilde{\xi}_i)_{i\in I}$ fulfills conjugate relations. Consider $n\geq 0$ and $i,i_1,\ldots, i_n\in I$ \st $(\tilde{\xi}_i, a_{i_1}, \ldots, a_{i_n})$ is a square sequence. Then \begin{eqnarray*}c_{n+1}(\tilde{\xi}_i\otimes a_{i_1}\otimes \cdots\otimes a_{i_n})&=&\f{\rho_{c(a_i)}}{\rho_{r(a_i)}} c_{n+1}(\xi_{*(i)}^{*}\otimes a_{*(i_1)}^*\otimes \cdots\otimes a_{*(i_n)}^*)\\ 
&=&\f{\rho_{c(a_i)}}{\rho_{r(a_i)}} (c_{n+1}(a_{*(i_n)}\otimes \cdots\otimes a_{ *(i_1)}  \otimes \xi_{*(i) })^*\\
 &=&(\eta_{r(\xi_i),c(\xi_i)}\circ c_{n+1}( \xi_{*(i) }\otimes a_{*(i_n)}\otimes \cdots\otimes a_{ *(i_1)} ))^*\\
 &=&(\eta_{r(\xi_i),c(\xi_i)}\delta_{n,1}\delta_{*(i),*(i_1)}p_{ r( \xi_{*(i)} ) })^*\\ &=&\delta_{i,i_1} p_{ r( \xi_{i} ) }
,\end{eqnarray*}
which closes the proof.
\end{pr}

3. A link with the already defined notions of Fisher's information can be made as follows: Consider $x\in \A_{kl}$ \st $x^*x$ is invertible in $\A_{ll}$. Then   a pair  $(\xi, \f{\rho_{k}}{\rho_{l}}\xi^*)$ fulfills conjugate relations for  $(x,x^*)$ if   
$\f{\rho_k}{\rho_l}\xi x(x^*x)^{-1}$ fulfills conjugate relations in the sens of \cite{NSS99} for $x^*x$ in $(\A_{ll},\vfi_l)$. The reciprocal is true if we have moreover $\vfi_l(\xi x(x^*x)^{-1})=0$. 

\begin{pr}$(\xi, \f{\rho_{k}}{\rho_{l}}\xi^*)$ fulfills conjugate relations for  $(x,x^*)$ \ssi for all $n\geq 0$, $$\ds\lf\{\begin{array}{crcl}(A)&\vfi_l(\xi x(x^*x)^n)&=\ds &\sum_{i=0}^n\vfi_k((xx^*)^i)\vfi_k((x^*x)^{n-i}),\\
(B)&\f{\rho_{k}}{\rho_{l}}\vfi_k(\xi^* x^*(xx^*)^n)&=\ds &\sum_{i=0}^n\vfi_l((x^*x^i)\vfi_k((xx^*)^{n-i}),\end{array}\ri.$$ where $(xx^*)^0$ stands for $p_k$ and  $(x^*x)^0$ stands for $p_l$.

But using $\vfi(.^*)=\overline{\vfi(.)}$, we have \begin{eqnarray*}(A)&\Leftrightarrow & \vfi_l((x^*x)^nx^*\xi^*)=\ds \sum_{i=0}^n\vfi_l((x^*x^i)\vfi_k((xx^*)^{n-i}),\end{eqnarray*}which is equivalent to $(B)$ because $\vfi_l((x^*x)^nx^*\xi^*)=\f{\rho_{k}}{\rho_{l}}\vfi_k(\xi^* x^*(xx^*)^n)$. 

So  $(\xi, \f{\rho_{k}}{\rho_{l}}\xi^*)$ fulfills conjugate relations for  $(x,x^*)$ \ssi for all $n\geq 0$, $$\ds\vfi_l(\xi x(x^*x)^n) =\sum_{i=0}^n\vfi_k((xx^*)^i)\vfi_k((x^*x)^{n-i}). $$ It is implied by the fact that  $\f{\rho_k}{\rho_l}\xi x(x^*x)^{-1}$ fulfills conjugate relations in the sens of \cite{NSS99} for $x^*x$ in $(\A_{ll},\vfi_l)$, and the reciprocal is true if  we have moreover $\vfi_l(\xi x(x^*x)^{-1})=0$. 
\end{pr} 

\subsection{Cram\'er-Rao inequality}
\begin{Th}Consider a non null element $a$ of $\A_{kl}$, with $k\neq l$ and $\rho_k\leq \rho_l$. Then $$\vfi(aa^*)\Phi_r(a,a^*)\geq \rho_k^2+\rho_l^2,$$ with equality \ssi there exists $c$ positive number \st the moments of $c.aa^*$ in $(\A_{kk}, \vfi_k)$ are the moments of the Marchenko-pastur distribution with parameter $\f{\rho_l}{\rho_k}$.
\end{Th} 

\begin{pr} If $(a, a^*)$ hos no conjugate system, then it is obvious. If a conjugate system for $(a,a^*)$ is $(\xi, \f{ \rho_{k} }{ \rho_{l} } \xi^*)$ (indeed, by (\ref{dudududadadaisallIwanttosaytoyou}), any conjugate system has this form), then \begin{eqnarray*}\Phi_r(a,a^*)&=&\lf(1+\f{\rho_{k}^2 }{\rho_{l}^2 }\ri) \vfi(\xi\xi^*),\end{eqnarray*}
so it suffices to prove that $$\vfi(aa^*)\vfi(\xi\xi^*)\geq\rho_l^2 $$
Note that we have $\vfi((\xi^*)^* a)=\rho_l$, so the result follows from the Cauchy-Schwarz inequality. 
Moreover, we have equality \ssi there exists $c>0$ \st $\xi^*=ca$, which is equivalent to the fact that for all $n\geq 2$ even, $$ c_n^{(l)}(a^*\otimes a\otimes\cdots \otimes a)=\delta_{n,2}\ff{c}$$ 
and  $$c_n^{(k)}(a\otimes a^*\otimes\cdots \otimes a^*)=\delta_{n,2} \f{ \rho_{l} }{ \rho_{k} }\ff{c},$$ which is equivalent,  by \teo \ref{regard.qui.perce.et.qui.met.quelque.chose.en.feu.dans.moi}, to the fact that the moments of $c.aa^*$ in $(\A_{kk}, \vfi_k)$ are the moments of the Marchenko-pastur distribution with parameter $\f{\rho_l}{\rho_k}$.
\end{pr}

\subsection{Fisher's information and freeness}
  \begin{Th} Consider a non null elements $x,y$ of $\A_{kl}$, with $k\neq l$. Then we have \begin{equation}\label{1.3.05.1} \Phi_r(x,y,x^*,y^*)\geq \Phi_r(x,x^*)+\Phi_r(y,y^*),\end{equation} and we have equality if $x,y$ are free with amalgamation over $\D$. Moreover, if $$ \Phi_r(x,y,x^*,y^*)= \Phi_r(x,x^*)+\Phi_r(y,y^*)<\infty,$$then $x,y$ are free with amalgamation over $\D$.  
\end{Th}

\begin{pr} - Let us prove that  $$\Phi_r(x,y,x^*,y^*)\geq \Phi_r(x,x^*)+\Phi_r(y,y^*).$$   If $\Phi_r(x,y,x^*,y^*)=\infty$, it is clear, and in the other case, let $(\xi, \f{\rho_k}{\rho_l}\xi^*, \zeta, \f{\rho_k}{\rho_l}\zeta^*)$ be the conjugate system for $(x,y,x^*,y^*)$. Then $(\xi, \f{\rho_k}{\rho_l}\xi^*)$, (resp. $(\zeta, \f{\rho_k}{\rho_l}\zeta^*)$) satisfy conjugate relations for $(x,x^*)$  (resp. for $(y,y^*)$), so the result is proved. 

- Suppose that $x,y$ are free with amalgamation over $\D$. If $\Phi(x,x^*)=\infty$ or $\Phi(y,y^*)=\infty$, then we have equality in (\ref{1.3.05.1}). In the other case, let  $(\xi, \f{\rho_k}{\rho_l}\xi^*)$, (resp. $(\zeta, \f{\rho_k}{\rho_l}\zeta^*)$) be a conjugate system for $(x,x^*)$  (resp. for $(y,y^*)$). It suffices to prove that $(\xi, \f{\rho_k}{\rho_l}\xi^*, \zeta, \f{\rho_k}{\rho_l}\zeta^*)$ is a conjugate system for $(x,y,x^*,y^*)$. It is clear that $$\xi, \f{\rho_k}{\rho_l}\xi^*, \zeta, \f{\rho_k}{\rho_l}\zeta^*\in \overline{\Alg (\{x,y,x^*,y^*\}\cup \D)}^{||.||_2},$$ so it suffices to prove that  $(\xi, \f{\rho_k}{\rho_l}\xi^*, \zeta, \f{\rho_k}{\rho_l}\zeta^*)$ fulfills conjugate relations for $(x,y,x^*,y^*)$.
Let us prove condition (iii) of definition \ref{Chao+Wozniak=Choose.the.Bosniaque}. 
Since  $(\xi, \f{\rho_k}{\rho_l}\xi^*)$, (resp. $(\zeta, \f{\rho_k}{\rho_l}\zeta^*)$) is a conjugate system for $(x,x^*)$  (resp. for $(y,y^*)$), we have $c_1(\xi)=c_1(\xi^*)=c_1(\zeta)=c_1(\zeta^*)=0$, and $c_2(\xi\otimes x)=c_2(\zeta\otimes y)=p_l$, $c_2(\f{\rho_k}{\rho_l}\zeta^*\otimes y^*)=c_2(\f{\rho_k}{\rho_l}\xi^*\otimes x^*)=p_k$. Since $\xi,\xi^*\in \overline{\Alg (\{x,x^*\}\cup \D)}^{||.||_2},$ and $\zeta,\zeta^*\in \overline{\Alg (\{y,y^*\}\cup \D)}^{||.||_2},$ 
by freeness with amalgamation over $\D$, we have $$c_2(\zeta\otimes x)=c_2(\zeta^*\otimes x^*)=c_2(\xi\otimes y)=c_2(\xi^*\otimes y^*)=0.$$
Consider now $n\geq 2$, and a square sequence $(T, a_1, \ldots, a_n)\in \{\xi, \f{\rho_k}{\rho_l}\xi^*, \zeta, \f{\rho_k}{\rho_l}\zeta^*\}\tii(\{x,x^*,y,y^*\}^n)$. Let us prove that  \begin{equation}\label{1.3.05.2}c_{n+1}(T\otimes a_1\otimes \cdots \otimes a_n)=0.\end{equation}  For example, we can suppose that $T=\xi$. If one of the $a_i$'s is $y$ or $y^*$, then (\ref{1.3.05.2}) is due to the freeness with amalgamation over $\D$.  If none of the $a_i$'s is $y$ or $y^*$, then (\ref{1.3.05.2}) is due to the fact that $(\xi, \f{\rho_k}{\rho_l}\xi^*)$ fulfills conjugate relations for $(x,x^*)$.

- In order to finish the proof, let us prove that if $(\xi, \f{\rho_k}{\rho_l}\xi^*)$ is a conjugate system for $(x,x^*)$, then we have  $x\xi=\xi^*x^*$. Both belong to $$p_k\overline{\Alg (\{x,x^*\}\cup\D)}^{||.||_2}p_k,$$ which is equal to 
$$  \overline{\{p_k\}\cup\{(xx^*)^n\ste n\geq 1\}}^{||.||_2}.$$Thus, since $\vfi_k$ is a faithfull trace state on $\A_{kk}$, it suffices that for all $n\geq 0$, $$\vfi_k((xx^*)^nx\xi)=\vfi_k ((xx^*)^n\xi^*x^*) ,$$ where  $(xx^*)^0$ stands for $p_k$. We have \begin{eqnarray*}\vfi_k((xx^*)^nx\xi)&=& \f{\rho_l}{\rho_k}\vfi_l(\xi(xx^*)^nx)\\
&=&\f{\rho_l}{\rho_k}\ds\sum_{i=0}^n\vfi_k((xx^*)^i)\vfi_l((xx^*)^{n-i}),
\end{eqnarray*}
which is a real number. Thus $$\vfi_k((xx^*)^n\xi^*x^*)=\vfi_k(\xi^*x^*(xx^*)^n)=\vfi_k(((xx^*)^nx\xi)^*)= \overline{ \vfi_k((xx^*)^nx\xi)}=\vfi_k((xx^*)^nx\xi).$$So we have proved that $x\xi=\xi^*x^*$.

- Now, suppose that $$ \Phi_r(x,y,x^*,y^*)= \Phi_r(x,x^*)+\Phi_r(y,y^*)<\infty.$$ By \teo \ref{article.tombant.en.ruine}, in order to prove the freeness with amalgamation over $\D$ of $x$ and $y$, it suffices to prove that for all $n\geq 2$, for all $z_1,\ldots, z_n$ taken in the algebras $\Alg(\{x,x^*\}\cup\D)$ and  $\Alg(\{y,y^*\}\cup\D)$, but not all in the same one, we have \begin{equation}\label{1.3.05.3}c_n(z_1\otimes\cdots\otimes z_n)=0.\end{equation}
By the formula of cumulants with products as entries (Theorem 2 of \cite{SS01}), 
we can suppose that $$z_1,\ldots, z_n\in \{x,x^*,y,y^*\}.$$ If the sequence $(z_1,\ldots, z_n)$ is not square, then (\ref{1.3.05.3}) holds by paragraph \ref{13.7.4.1} (a). So we suppose the sequence to be square, and by equation (\ref{ecologie&feminisme}), we can suppose that $$(z_1,z_2)\in \{(x,y^*), (y,x^*), (x^*, y), (y^*,x)\}.$$For example, we will treat the case where $(z_1,z_2)=(x,y^*)$.

Consider  $(\xi, \f{\rho_k}{\rho_l}\xi^*, \zeta, \f{\rho_k}{\rho_l}\zeta^*)$ a conjugate system for $(x,y,x^*,y^*)$. Then using the  hypothesis and 1. of   \ref{retard.1.3.05.1}, we know that $(\xi, \f{\rho_k}{\rho_l}\xi^*)$, (resp. $(\zeta, \f{\rho_k}{\rho_l}\zeta^*)$) is the conjugate system for $(x,x^*)$  (resp. for $(y,y^*)$). 
Since  $x\xi=\xi^*x^*$, we have \begin{equation}\label{bleu.1.3.05}c_{n+1}^{(k)}(x\xi\otimes z_2\otimes\cdots \otimes z_n\otimes z_1^*)=c_{n+1}^{(k)}(\xi^*x^*\otimes z_2\otimes\cdots \otimes z_n\otimes z_1^*). \end{equation}
Now, we apply the   formula of cumulants with products as entries (Theorem 2 of \cite{SS01}) to left hand side and right hand side of (\ref{bleu.1.3.05}). 
 \begin{eqnarray*}LHS&=& c_{n+2}^{(k)}(x\otimes \xi \otimes y \otimes z_3\otimes \cdots \otimes z_n\otimes z_1^*)+\sum_{\substack{1\leq i\leq n \\ \textrm{$i$ even}}}\\ &&\underbrace{c_i^{(l)}(\xi\otimes z_2 \otimes z_3\otimes \cdots \otimes z_i)}_{\substack{=0\textrm{ because}\\ z_2=y\textrm{ (and $c_2(\xi)=0$)}}}c_{n-i+2}^{(k)}(x\otimes z_{i+1}\otimes \cdots \otimes z_n\otimes z_1^*)\\ &=& \f{\rho_l}{\rho_k}c_{n+2}^{(l)}(\xi \otimes y \otimes z_3\otimes \cdots \otimes z_n\otimes z_1^*\otimes x) \\ &=& 0.\end{eqnarray*}On the other side, we have 
\begin{eqnarray*}
RHS&=& \underbrace{c_{n+2}^{(k)}(\xi^*\otimes x^* \otimes z_2 \otimes z_3\otimes \cdots \otimes z_n\otimes z_1^*)}_{\substack{=0\textrm{ because}\\ n+2>2}}
+\sum_{\substack{1\leq i\leq n \\ \textrm{$i$ even}}}\\ &&c_i^{(l)}(x^*\otimes z_2 \otimes z_3\otimes \cdots \otimes z_i)\underbrace{c_{n-i+2}^{(k)}(\xi^*\otimes z_{i+1}\otimes \cdots \otimes z_n\otimes z_1^*)}_{=0\textrm{ if $i<n$}}\\ &=& c_n^{(l)}(x^*\otimes y\otimes z_3\otimes \cdots \otimes z_n)\underbrace{c_2^{(k)}(\xi^*\otimes x^*)}_{ \neq  0}.\end{eqnarray*}Hence $$c_n^{(l)}(x^*\otimes y\otimes z_3\otimes \cdots \otimes z_n)=0,$$ which is what we wanted to obtain.
\end{pr}
 
\subsubsection*{Question} It would be interesting to have a characterization of $R$-diagonal elements with non trivial kernel involving the Fisher's information defined in this paper. Inspired by the paper \cite{NSS99}, we ask the following question: given a compactly supported \pro measure $\nu$ on $\R^+$, what are the elements $a\in \A p_1$ \st $aa^*$ has distribution $\nu$ in $(\A_{22},p_2)$, and \st  $\Phi_r(p_1ap_{2},p_2ap_{2}, p_2a^*p_{1}, p_2a^*p_{2})$ is minimal ?

\section*{Appendix A: proof of proposition \ref{livre.cardio.pauline}}
Let us suppose, for example, that $\rho_k \leq \rho_l$.  Fix $R>0$ \st $R^2$ is more than the suppremum of the support of $\mu$. 
Define, for all $q\geq 1$, the map $$\kappa_q\, :\, x\in [0,R]^q\to \ff{q}\ds\sum_{i=1}^q\delta_{x_i}.$$ The index $q$ in $\kappa_q$ will always be omitted, because no confusion will ever be possible. For $P$ \pro measure on $[0,R]^q$, we denote by $\kappa(P)$ the push-forward of $P$ by $\kappa$.  
First of all, let us recall a large deviation principle. For basic definitions on large deviations, see P. 177 of \cite{hiai} or many other books (e.g. \cite{DemboZeitouni}). 

\begin{Th}Let, for $n\geq 1$, $Z_n$ be the total mass of 
$$P_n:= 
\ds\Delta(x)^2  
\prod_{i=1}^{q_k(n)}x_j^{q_l(n)-q_k(n)} 
1_{[0,R]^{q_k(n)}}(x)\ud x.$$  Then the finite limit $B:=\ds\lim_{n\to\infty}n^{-2}\log Z_n$ exists, and the sequence  $\lf(\kappa\lf(\ff{Z_n}P_n\ri)\ri)_n$ satisfies a large deviation principle in the set of \pro measures on $[0,R]$ endowed with topology of weak convergence in the scale $n^{-2}$ with the good rate function $$I\,:\,\nu\to -\rho_k^2\Sigma(\nu)-\lf(\rho_k\rho_l-\rho_k^2\ri)\int \log(x)\ud \nu(x)+B.$$  
\end{Th}
This 
\teo 
was 
proved under a slightly different hypothesis in \cite{hiai} (\teo 5.5.1 p. 227, with $Q=0$). The difference between the hypothesis above and the hypothesis of  \teo 5.5.1 of \cite{hiai} is that in the latter, the bound $R$ does not appear, the measures are considered on $\R^+$. But it is not a problem: the proof of \teo 5.5.1 can easily be adaptated to this context (in fact it is more easy to work with the compact set $[0,R]$). Note that an analoguous modification of a result proved for the interval $\R^+$ to the interval $[0,R]$  is done p. 240 of \cite{hiai}. 

Note that removing the renormalization constant $Z_n$ and the limit $B$, one gets the following result. Its formulation  implies to have extended the notion of large deviation principle to sequences of finite measures (not only of \pro measures), but it can be done without any ambiguity.

\begin{cor}The sequence  of finite measures $\lf(\kappa\lf(P_n\ri)\ri)_n$ satisfies a large deviation principle in the set of \pro measures on $[0,R]$ endowed with topology of weak convergence in the scale $n^{-2}$ with the good rate function $$J\,:\,\nu\to -\rho_k^2\Sigma(\nu)-\lf(\rho_k\rho_l-\rho_k^2\ri)\int \log(x)\ud \nu(x).$$  \end{cor}

Now, we give the proof of proposition \ref{livre.cardio.pauline}. 

{\it Step I. } For all $r\geq 1$, $\eps>0$, 
$$\chi_R(a;2r,\eps)=\ds\limsup_{n\to\infty} \ff{n^2}\log\La (\Gamma_R(a;n,2r,\rho_k\eps))+\rho_k \rho_l\log n .$$ 
$\Gamma_R(a;n,r,\rho_k\eps)$ 
is the set of matrices of $\Mn (k,l)$ \st $||M||\leq R$ and each moment of order $\leq r$ of the spectral law of the $k$-th diagonal block of $MM^*$ is $\eps$-closed to the  moment of same order of $\mu$. 
Thus by the remark following lemma \ref{16.2.05.2}, $\La(\Gamma_R(a;n,2r,\rho_k\eps))$
is $$\ds\f{\pi^{q_k(n)q_l(n)}}{\prod_{j=1}^{q_k(n)}j!\prod_{j=q_l(n)-q_k(n)}^{q_l(n)-1}j!} P_n(\{x\in [0,R]^{q_k(n)}\ste \forall i=0,\ldots, r, |m_i(\kappa(x))-m_i(\mu)|\leq \eps\}),$$where $P_n$ is the measure introduced in the previous theorem.

{\it Step II. } Let us compute the limit $C$, as $n\to \infty$, of $$\ds u_n:=\ff{n^2}\log \f{\pi^{q_k(n)q_l(n)}}{\prod_{j=1}^{q_k(n)}j!\prod_{j=q_l(n)-q_k(n)}^{q_l(n)-1}j!}+\rho_k\rho_l\log n.$$
We have, by Stirling formula, $\log j! = \ff{2}\log j +j(\log j-1) +O(1).$ So 
$$u_n= \f{q_k(n)q_l(n)}{n^2}\log \pi -\ff{n^2}\ds\sum_{j=1}^{q_k(n)}(\ff{2}\log j +j(\log j-1)) -\ff{n^2}\!\!\!\!\sum_{j=q_l(n)-q_k(n)}^{q_l(n)-1}\!\!\!\!(\ff{2}\log j +j(\log j-1))$$ $$+\rho_k\rho_l\log n +O(\ff{n}) $$ $$=
\rho_k\rho_l\log\pi+\f{q_k(n)(q_k(n)+1)+(q_l(n)-1)q_l(n)-(q_l(n)-q_k(n)-1)(q_l(n)-q_k(n))}{2n^2}$$ 
$$ 
-\f{q_k(n)^2}{n^2}\ds\ff{q_k(n)}\sum_{j=1}^{q_k(n)}\f{j}{q_k(n)}\log \f{j}{q_k(n)} -\f{\log q_k(n)}{n^2}\ds\sum_{j=1}^{q_k(n)}j $$ 
$$
-\f{ q_k(n)^2 }{ n^2 }\ds\ff{q_k(n)}\!\!\sum_{j=q_l(n)-q_k(n)}^{q_l(n)-1}\!\!\f{j}{q_k(n)}\log \f{j}{q_k(n)} -\f{\log q_k(n)}{n^2}\ds\!\!\!\!\sum_{j=q_l(n)-q_k(n)}^{q_l(n)-1}\!\!\!\! j  $$     $$+\rho_k\rho_l\log n  +o(1)
$$ $$=\rho_k\rho_l(\log\pi+1)-\rho_k^2\int_0^1t\log t\ud t -\rho_k^2\int_{\rho_l/\rho_k-1}^{\rho_l/\rho_k}t\log t\ud t $$ $$-\log q_k(n)\f{q_k(n)(q_k(n)+1)+(q_l(n)-1)q_l(n)-(q_l(n)-q_k(n)-1)(q_l(n)-q_k(n))}{2n^2}$$  $$+\rho_k\rho_l\log (q_k(n)\f{n}{q_k(n)})  +o(1)$$ 
$$=\rho_k\rho_l(\log\pi+1-\log \rho_k)-\rho_k^2\int_0^1t\log t\ud t -\rho_k^2\int_{\rho_l/\rho_k-1}^{\rho_l/\rho_k}t\log t\ud t +\log q_k(n)\tii $$ $$\lf( \f{q_k(n)(q_k(n)+1)+(q_l(n)-1)q_l(n)-(q_l(n)-q_k(n)-1)(q_l(n)-q_k(n))-2q_k(n)q_l(n)+O(n)}{2n^2}\ri)$$ $$+o(1).$$
Thus $$u_n\ninf C:=\rho_k\rho_l(\log\pi+1-\log \rho_k)-\rho_k^2\underbrace{\int_0^1t\log t\ud t}_{=-1/4} -\rho_k^2\int_{\rho_l/\rho_k-1}^{\rho_l/\rho_k}t\log t\ud t. $$

{\it Step III. } 
Now, let us denote by $F(r,\eps)$ (resp. $G(r,\eps)$) the set of \pro measures on $[0,R]$ for which each moment of order $\leq r$ is $\eps$-closed (resp. strictly $\eps$-closed) to the moment of same order of $\mu$.  $F(r,\eps)$ (resp. $G(r,\eps)$) is closed (resp. open). Thus by the previous corollary, we have $$\ds\limsup_{n\to\infty}\ff{n^2} \log P_n(\kappa^{-1}(F(r,\eps)))\leq -\underset{F(r,\eps)}{\inf} J,$$
$$\ds\liminf_{n\to\infty}\ff{n^2} \log P_n(\kappa^{-1}(G(r,\eps)))\geq -\underset{G(r,\eps)}{\inf} J.$$
But $\underset{F(r,\eps)}{\inf} J=\underset{G(r,\eps)}{\inf} J,$ so $$\ff{n^2} \log P_n(\kappa^{-1}(F(r,\eps)))\ninf -\underset{F(r,\eps)}{\inf} J,$$ and it follows, by Steps I and II,  that $$\chi_R(a;2r,\eps)=\ds -\underset{F(r,\eps)}{\inf} J+C.$$ As $\eps $ goes to $0$ and $r$ goes to $\infty$, $\ds\underset{F(r,\eps)}{\inf} J$ goes to $J(\mu)$, and we obtain the desired result: $$\chi_R(a)=\rho_k^2\Sigma(\mu) +(\rho_l-\rho_k)\rho_k\int\log(x)\ud\mu(x)+\rho_k\rho_l(\log\f{\pi}{\rho_k}+1)+\f{\rho_k^2}{4}-\rho_k^2\int_{\rho_l/\rho_k-1}^{\rho_l/\rho_k}t\log t\ud t. $$

\section*{Appendix B: proof of lemma \ref{16.2.05.2}} 
In all this proof, we shall identify elements of $\R^q$ with the associated diagonal $q\tii q$ matrix.

a) First of all, the fact that $\Psi$ is an injection onto a set a negligeable complementary is well known (see \cite{H&J}).

b) Let $P$ be the $q\tii q'$ matrix with entry $(i,j)$ equal to $1$ if $i=j$, and to $0$ in the other case. Then the set $\mc{U}_{q,q'}$ is $P\mc{U}_{q'}$. So the set 
 $\mc{U}_q/T_q\,\tii\, \R^q_{+,<}\,\tii\, \mc{U}_{q,q'}$ is a manifold, and for $u\in \mc{U}_q$, $ a\in\R^q_{+,<}$, $v\in \mc{U}_{q'}$,  its tangent space at $(uT_q,a, Pv)$ is the cartesian product of tangent spaces of  respectively $\mc{U}_q/T_q,\, \R^q_{+,<},\, \mc{U}_{q,q'}$ at respectively $uT_q,a, Pv$. The first of them can be identified, via the map $M\to u^*M$, to the set $\mf{U}^0_q$ of anti-hermitian matrices with zeros an the diagonal, the second one is $\R^n$, and the third one can be identified, via   the map $M\to Mv^*$,   to the set $\mf{U}_{q,q'}$ of $q\tii q'$ matrices  in which the submatrix of the first $q$ columns is anti-hermitian. The differential $\ud \Psi_{uT_q,a, Pv}$ of $\Psi $  at $(uT_q,a, Pv)$ is given by the following formulae: \begin{eqnarray*}\forall X\in \mf{U}^0_{q},\quad \ud \Psi_{uT_q,a, Pv}(u^*X,0,0)&=&u(\underbrace{Xa^{1/\!2}-a^{1/\!2}X}_{:=M_a(X)})u^*Pv,\\  \forall A\in \R^q,\quad  \ud \Psi_{uT_q,a, Pv}(0,A,0)&=&\ff{2}u\f{A}{a^{1/\!2}}u^*Pv,\\
\forall Y\in \mf{U}_{q,q'}, \quad \ud \Psi_{uT_q,a, Pv}(0,0,Yv)&=&ua^{1/\!2}u^*Yv\end{eqnarray*}

c) Let $\det$ be the determinant on the canonical basis of $\mf{M}_{q,q'}$, i.e. on the basis containing the  elementary matrices and $i$ times this matrices. Define $n=q^2-q$ and $m=2qq'-q^2$. Fix  $(uT_q,a, Pv)$ in the manifold, and $X_1$,\ldots, $X_n\in \mf{U}^0_{q}$, $A_1,$\ldots ,$A_q\in \R^q$, and $Y_1$, \ldots, $Y_m\in \mf{U}_{q,q'}$. Now, let us compute the differential form $\Psi^*\!\det$ at $(uT_q,a, Pv)$ on the family $$((u^*X_1,0,0),\ldots, (u^*X_n,0,0),(0,A_1,0),\ldots ,(0,A_q,0), (0,0,Y_1v), \ldots, (0,0,Y_mv)).$$ It is $$\det(uM_a(X_1)u^*Pv,\ldots, uM_a(X_n)u^*Pv,\ff{2}u\f{A_1}{a^{1/\!2}}u^*Pv,\ldots ,\ff{2}u\f{A_q}{a^{1/\!2}}u^*Pv,ua^{1/\!2}u^*Y_1v , \ldots, ua^{1/\!2}u^*Y_mv).$$
Define $\tilde{u}=\begin{bmatrix}u&0\\ 0&I_{q'-q}\end{bmatrix}\in \mc{U}_{q'}$. We have $P\tilde{u}=uP$.  Note that the base we choosed is orthonormal for the euclidian structure we choosed on $\mf{M}_{q,q'}$, and the left or right multiplications by unitary elements are orthogonal, and have determinant $1$ by connexity of the unitary group. So  what we want to compute is  equal to$$\det (M_a(X_1)P,\ldots, M_a(X_n)P,\ff{2}\f{A_1}{a^{1/\!2}}P,\ldots ,\ff{2}\f{A_q}{a^{1/\!2}}P,a^{1/\!2}u^*Y_1\tilde{u} , \ldots, a^{1/\!2}u^*Y_m\tilde{u}).$$
In order to compute this, let us introduce another basis of $\mf{M}_{q,q'}$. Let $(E_{k,l})_{\substack{1\leq k\leq q\\ 1\leq l\leq q'}}$
 be the elementary matrices of $\mf{M}_{q,q'}$. Define, for $1\leq k<l\leq q$, 
$$e_{k,l}=E_{k,l}+E_{l,k},\quad e'_{k,l}=i(E_{k,l}-E_{l,k}),$$let $\mc{B}_1$ be the family $(e_{k,l}, e'_{k,l})_{1\leq k<l\leq q }$, define $$\mc{B}_2=(E_{k,k})_{1\leq k\leq q},$$  and let $\mc{B}_3$ be any basis of  $\mf{U}_{q,q'}$. Note that $\mc{B}:= \mc{B}_1\cup\mc{B}_2\cup\mc{B}_3$ is a basis of $\mf{M}_{q,q'}$. Let $\la$ be a non-null real number \st $\det$ is $\la$ times the determinant on $\mc{B}$. Note that the matrix of the family $$\mc{F}:= (M_a(X_1)P,\ldots, M_a(X_n)P,\ff{2}\f{A_1}{a^{1/\!2}}P,\ldots ,\ff{2}\f{A_q}{a^{1/\!2}}P,a^{1/\!2}u^*Y_1\tilde{u} , \ldots, a^{1/\!2}u^*Y_m\tilde{u})$$on $\mc{B}$ is block upper-triangular (with respect to the decomposition $\mc{B}= \mc{B}_1\cup\mc{B}_2\cup\mc{B}_3$). So $\det(\mc{F})$ is $\la$ times the product of the determinants of the matrices of the families $$ (M_a(X_1)P,\ldots, M_a(X_n)P),\;\; \;(\ff{2}\f{A_1}{a^{1/\!2}}P,\ldots ,\ff{2}\f{A_q}{a^{1/\!2}}P),\; \; \;(\Proj(a^{1/\!2}u^*Y_1\tilde{u}) , \ldots, \Proj(a^{1/\!2}u^*Y_m\tilde{u}))$$  on the respective bases $\mc{B}_1$,  $\mc{B}_2$, $\mc{B}_3$ (where $\Proj$ denotes the projection on $\Span(\mc{B}_3)=\mf{U}_{q,q'}$ in the direction of $\Span(\mc{B}_1\cup\mc{B}_2)$, i.e. the orthogonal projection on $\mf{U}_{q,q'}$). 

Let us compute the first determinant. $M_a$ maps linearly  $\mf{U}^0_{q}$ into the set of hermitian matrices with null diagonal, so $X\to M_a(X)P$ maps linearly  $\mf{U}^0_{q}$ into $\Span(\mc{B}_1)$. Let   $(F_{k,l})_{\substack{1\leq k,l\leq q}}$
 be the elementary matrices of $\mf{M}_{q,q}$. Define, for $1\leq k<l\leq q$, 
$$f_{k,l}=F_{k,l}-F_{l,k},\quad f_{k,l}'=i(F_{k,l}+F_{l,k}),$$let $\beta_1$ be the family $(f_{k,l}, f'_{k,l})_{1\leq k<l\leq q }$. The matrix of the map $X\to M_a(X)P$ between the bases  $\beta_1$ and  $\mc{B}_1$ is block-diagonal, with blocks $$\begin{array}{cr}\begin{bmatrix}0&a^{1/\!2}_k-a^{1/\!2}_l\\ a^{1/\!2}_l-a^{1/\!2}_k&0\end{bmatrix},&\quad\quad (1\leq k<l\leq q).\end{array}$$
So its determinant is $\Delta(a^{1/\!2})^2$, and the determinant of the matrix  of the family $ (M_a(X_1)P$,\ldots, $M_a(X_n)P)$ in 
$\mc{B}_1$ is $\Delta(a^{1/\!2})^2$ times the determinant of the matrix  of the family $ (X_1,\ldots, X_n)$ in 
$\beta_1$. 

 The second determinant 
is $\ff{2^n}\ff{(a_1\cdots a_n)^{1/\! 2}}$ times  the determinant of the matrix  of the family $ (A_1,\ldots, A_n)$ in 
$\mc{B}_2$.  

Let us compute the third determinant. In this paragraph, we shall use a two blocks-decomposition of matrices of $\mf{M}_{q,q'}$. Any $q\tii q'$ matrix $Y$ will be denoted by $Y=(Y_s,Y_r)$, where $Y_s$ is a $q\tii q$ matrix, and $Y_r$ is a $q\tii (q'-q)$ matrix ($s$ is for {\it square}, and $r$ for {\it rectangular}). With this decomposition, $\Proj$ has a simple expression: $\Proj(Y)=(\f{Y_s-Y_s^*}{2},Y_r)$. Note that if $Y\in \mf{U}_{q,q'}$, then $Y_s$ is anti-hermitian, so $$\Proj(a^{1/\!2}u^*Y\tilde{u})=N_a((u^*Y_su,u^*Y_r)),$$where $$N_a \, :\, Y\in \mf{U}_{q,q'} \to (\f{a^{1/\!2}Y_s+Y_sa^{1/\!2}}{2},a^{1/\!2}Y_r)\in \mf{U}_{q,q'}.$$ Note that $Y\in  \mf{U}_{q,q'}\to (u^*Y_su,u^*Y_r)$ is orthogonal, and has determinant one by connexity of $\mc{U}_q$. Let us compute the determinant  of $N_a$.  All vectors of the basis $$(E_{k,l}-E_{l,k})_{1\leq k<l\leq q}\cup      (i(E_{k,l}+E_{l,k}))_{1\leq k\leq l\leq q}\cup    (E_{k,l})_{\substack{\!\!1\leq k\leq q\\ q'-q\leq l\leq q'}}\cup      (iE_{k,l})_{\substack{\!\!1\leq k\leq q\\ q'-q\leq l\leq q'}}   $$ are eigenvectors of $N_a$, with respective eigenvalues
  $$(\f{a_k^{1/\! 2}+a_l^{1/\! 2}}{2})_{1\leq k<l\leq q}\cup      (\f{a_k^{1/\! 2}+a_l^{1/\! 2}}{2})_{1\leq k\leq l\leq q}\cup    (a_k^{1/\! 2})_{\substack{1\leq k\leq q\\ q'-q\leq l\leq q'}}\cup      (a_k^{1/\! 2})_{\substack{1\leq k\leq q\\ q'-q\leq l\leq q'}}   .$$ 
Thus the determinant of $N_a\, :\,\mf{U}_{q,q'} \to\mf{U}_{q,q'}$ is$$(a_1\cdots a_q)^{q'-q+\ff{2}}\ds\prod_{1\leq k<l\leq q}\lf(\f{a_k^{1/\! 2}+a_l^{1/\! 2}}{2}\ri)^2.$$The   third determinant is this quantity times the determinant of the family  $(Y_1, \ldots, Y_m)$ in $\mc{B}_3$.

d) So the value of the differential form $\Psi^*\!\det$ at $(uT_q,a, Pv)$ on the family $$((u^*X_1,0,0),\ldots, (u^*X_n,0,0),(0,A_1,0),\ldots ,(0,A_q,0), (0,0,Y_1v), \ldots, (0,0,Y_mv))$$ 
is $$\la \Delta(a^{1/\!2})^2\det_{\beta_1}(X_1,\ldots, X_n)\ff{2^n}\ff{(a_1\cdots a_n)^{1/\! 2}}\tii $$ $$\det_{\textrm{can}}(A_1,\ldots, A_n)(a_1\cdots a_q)^{q'-q+\ff{2}}\ds\prod_{1\leq k<l\leq q}\lf(\f{a_k^{1/\! 2}+a_l^{1/\! 2}}{2}\ri)^2\det_{\mc{B}_3}(Y_1, \ldots, Y_m).$$ It is well known (see, e.g., section I.5 of \cite{BrockerDieck}), that it  is equal, up to a multiplicative constant, to $$\Delta(a)^2\ds\prod_{k=1}^qa_k^{q'-q}\omega^{\mc{U}_q/T_q}_{uT}(u^*X_1,\ldots, u^*X_n)\det_{\textrm{can}}(A_1,\ldots, A_n)\omega^{\mc{U}_{q,q'}}_{vP}(Y_1v,\ldots, Y_mv),$$ where $\omega^{\mc{U}_q/T_q}$ is a non-null differential $n$-form on $\mc{U}_q/T_q$ which is invariant under the left action of the unitary group, and  $\omega^{\mc{U}_{q,q'}}$ is a non-null differential $m$-form on $\mc{U}_{q,q'}$ which is invariant under the left and right actions of the unitary groups. So $\Psi^*\!\det$ is equal, up to a multiplicative constant, to $$f.\omega^{\mc{U}_q/T_q}\wedge \det_{\textrm{can}}\wedge\,\omega^{\mc{U}_{q,q'}}, $$ where $f$ is the smooth function defined on $\mc{U}_q/T_q\,\tii\, \R^q_{+,<}\,\tii\, \mc{U}_{q,q'}$ by $$f(uT, a, Pv)=\Delta(a)^2\ds\prod_{k=1}^qa_k^{q'-q}.$$ Hence
the push-forward, by $\Psi^{-1}$, of the Lebesgue measure on $\mf{M}_{q,q'}$ is the tensor product $\gamma_q\otimes C\sigma_{q,q'}\otimes \gamma_{q,q'},$ where $C$ is a positive constant.

e) Let us conpute $C$. As noticed in the remark following the lemma, by definition of the measures $\gamma_q$ and $\gamma_{q,q'}$, we can now claim that the map \begin{eqnarray*}\tilde{\Psi}\, :\,\mc{U}_q\,\tii\, (\R^+)^q\,\tii\, \mc{U}_{q'}&\to& \mf{M}_{q,q'} \\
(u,x,v)&\mapsto & u\diag (x_1,\ldots, x_q)^{1/\! 2}u^*Pv\end{eqnarray*}
is surjective and preserves the measure $\Haar (\mc{U}_q)\otimes \f{C}{q!}\tilde{\sigma}_{q,q'}\otimes \Haar (\mc{U}_{q'})$ (i.e. the push-forward of this measure by $\tilde{\Psi}$ is the Lebesgue measure), where $\tilde{\sigma}_{q,q'}$ is the measure on $(\R^+)^n$ with density given by formula (\ref{banque.avec.david}).  

The function $ x\in \mf{M}_{q,q'}\to e^{-\Tr xx^*}$ has integral with respect to the Lebesgue measure equal to $\pi^{qq'}$. Thus
\begin{eqnarray*}\pi^{qq'}&=&  \ds\f{C\pi^{qq'}}{\prod_{j=1}^{q}j!\prod_{j=q'-q}^{q'-1}j!}\int_{a\in (\R^+)^q}\int_{u\in \mc{U}_q}\int_{v\in \mc{U}_{q'}}\Delta(a)^2\prod_{j=1}^{q}a_j^{q'-q}e^{-\Tr ua^{1/\! 2}Pvv^*P^*a^{1/\! 2}u^*}\ud a\ud u\ud v. \end{eqnarray*}
Thus \begin{eqnarray*}\ff{C}&=&  \ds\ff{\prod_{j=1}^{q}j!\prod_{j=q'-q}^{q'-1}j!}\int_{a\in (\R^+)^q}\Delta(a)^2\prod_{j=1}^{q}a_j^{q'-q}e^{-\sum_{i=1}^qa_i}\ud a \end{eqnarray*}
We can now apply formula (4.1.8) p. 119 of \cite{hiai}, with $n=q$, $\beta=1$, $a=q'-q+1$, and it appears that $C=1$.

Florent Benaych-Georges\\ DMA, \'Ecole Normale Sup\'erieure,\\ 45 rue d'Ulm, 75230 Paris Cedex 05, France\\   e-mail: benaych@dma.ens.fr\\
  http://www.dma.ens.fr/$\sim$benaych
 \end{document}